\documentclass[article,10pt,pdflatex,reqno]{amsart}
\usepackage{amsmath,times}
\usepackage{amsthm}
\usepackage{amssymb}
\usepackage{bbm}
\usepackage{latexsym}
\usepackage{amscd}
\usepackage[english]{babel}
\usepackage[latin1]{inputenc}
\usepackage{graphicx}
\usepackage{color}
\usepackage{verbatim}
\newtheorem{theorem}{Theorem}[section]
\newtheorem{lemma}[theorem]{Lemma}
\newtheorem{proposition}[theorem]{Proposition}

\newtheorem*{thmA}{Theorem A} 
\newtheorem*{thmAp}{Theorem A'} 
\newtheorem*{thmB}{Theorem B} 
\newtheorem*{thmBp}{Theorem B'}
\newtheorem*{thmC}{Theorem C}

\newtheorem*{corA}{Corollary A}

\newtheorem*{ackn}{Acknowledgment} 

\theoremstyle{definition}
\newtheorem{definition}[theorem]{Definition}
\newtheorem{corollary}[theorem]{Corollary}
\newtheorem{remark}[theorem]{Remark}

\begin{document} 

\title{Deformations of K\"ahler manifolds to normal bundles and restricted volumes of big classes}
\author{David Witt Nystr\"om}
\maketitle

\begin{abstract}
The deformation of a variety $X$ to the normal cone of a subvariety $Y$ is a classical construction in algebraic geometry. In this paper we study the case when $(X,\omega)$ is a compact K\"ahler manifold and $Y$ is a submanifold. The deformation space $\mathcal{X}$ is fibered over $\mathbb{P}^1$ and all the fibers $X_{\tau}$ are isomorphic to $X$, except the zero-fiber, which has the projective completion of the normal bundle $N_{Y|X}$ as one of its components. The first main result of this paper is that one can find K\"ahler forms on modifications of $\mathcal{X}$ which restricts to $\omega$ on $X_1$ and which makes the volume of the normal bundle in the zero-fiber come arbitrarily close to the volume of $X$. Phrased differently, we find K\"ahler deformations of $(X,\omega)$ such that almost all of the mass ends up in the normal bundle. The proof relies on a general result on the volume of big cohomology classes, which is the other main result of the paper. A $(1,1)$ cohomology class on a compact K\"ahler manifold $X$ is said to be big if it contains the sum of a K\"ahler form and a closed positive current. A quantative measure of bigness is provided by the volume function, and there is also a related notion of restricted volume along a submanifold. We prove that if $Y$ is a smooth hypersurface which intersects the K\"ahler locus of a big class $\alpha$ then up to a dimensional constant, the restricted volume of $\alpha$ along $Y$ is equal to the derivative of the volume at $\alpha$ in the direction of the cohomology class of $Y$. This generalizes the corresponding result on the volume of line bundles due to Boucksom-Favre-Jonsson and independently Lazarsfeld-Musta\c{t}\u{a}.   
\end{abstract}

\section{Introduction}

\subsection{Deformation to the normal cone}

Let $X$ be a complex manifold (or variety) and $Y$ a subvariety of $X$. The blow-up $\mu: \mathcal{X}\to X\times \mathbb{P}^1$ of $Y\times \{0\}$ is a classical construction in algebraic geometry known as the \emph{deformation to the normal cone of $Y$} (see e.g. \cite{Ful}). The zero fiber $X_0$ has two components: $X'$ which is isomorphic to the blow-up $\mu'$ of $Y$ in $X$, and the exceptional divisor $\mathcal{E}$ which is  naturally identified with the projective completion $\mathbb{P}(N_{Y|X}\oplus \mathbb{C})$ of the normal cone $N_{Y|X}$ of $Y$. The intersection $\mathcal{E}\cap X'$ identifies the exceptional divisor $E$ in $X'$ with the divisor at infinity of $\mathbb{P}(N_{Y|X}\oplus \mathbb{C})$, so $X_0\setminus X'=N_{Y|X}$.

Since $Y\times \{0\}$ is invariant under the $\mathbb{C}^*$-action on $X\times \mathbb{P}^1$ which acts by multiplication on the base, this action lifts to $\mathcal{X}$. On $\mathcal{E}=\mathbb{P}(N_{Y|X}\oplus \mathbb{C})$ the action is that induced by multiplication on the factor $\mathbb{C}$, and hence on $N_{Y|X}$ it acts by inverse multiplication.

\subsection{Deforming a K\"ahler manifold}

We are interested in the case when $(X,\omega)$ is a compact K\"ahler manifold, and $Y$ is a submanifold of $X$. $\mathcal{X}$ is then also compact K\"ahler, and $N_{Y|X}$ is simply the normal bundle of $Y$.

If $\Omega$ is an $S^1$-invariant K\"ahler form on $\mathcal{X}$ such that $\Omega_{|X_1}=\omega$ it encodes a deformation of $(X,\omega)$ to the singular space $(\mathcal{E}, \Omega_{|\mathcal{E}})\cup (X',\Omega_{|X'})$. We call $(\mathcal{X},\Omega)$ a \emph{K\"ahler deformation} of $(X,\omega)$ to the normal bundle of $Y$.

Let $n$ denote the complex dimension of $X$. Clearly $$\int_{N_{Y|X}}\Omega^n+\int_{X'}\Omega^n=\int_{\mathcal{E}}\Omega^n+\int_{X'}\Omega^n=\int_X \omega^n.$$  

Can we always find K\"ahler deformations $(\mathcal{X},\Omega)$ to the normal bundle of $Y$ so that almost all of the volume of $(X,\omega)$ ends up in $N_{Y|X}$?

The answer is no. To see why, consider the $(1,1)$ cohomology class $\beta:=[\Omega]$ associated to a deformation $(\mathcal{X},\Omega)$. The cohomology of $\mathcal{X}$ is generated by the pullback of the cohomology of $X$ together with the classes $\{X_0\}$ and $\{\mathcal{E}\}$ corresponding to the currents of integration $[X_0]$ and $[\mathcal{E}]$. If $\alpha:=[\omega]$ we have that $\beta_{|X_1}=\alpha$ and so $$\beta=(\pi_X\circ\mu)^*\alpha+b\{X_0\}-c\{\mathcal{E}\}$$ for some constants $b,c$. Consequently $\beta_{|X'}=\mu'^*\alpha-c\{E\}$. 

Recall that a $(1,1)$-class is called K\"ahler if it contains a K\"ahler form. Since $\Omega_{|X'}$ is K\"ahler the class $\mu'^*\alpha-c\{E\}$ must be K\"ahler. The supremum of $t$ such that $\mu'^*\alpha-t\{E\}$ is K\"ahler is known as the \emph{Seshadri constant} $\epsilon(\alpha,Y)$ of $Y$.

The function $\int_{X'}(\mu'^*\alpha-t\{E\})^n$ of $t\in (0,\epsilon(\alpha,Y))$ is decreasing and its infimum is equal to $\int_{X'}(\mu'^*\alpha-\epsilon(\alpha,Y)\{E\})^n$. The issue is that this number can be strictly positive. 

The easiest example of this happening is probably $X=\mathbb{P}^1\times \mathbb{P}^1$, $Y=\{0\}\times \{0\}$ and $\alpha=[\pi_1^*\omega_{FS}+\pi_2^*\omega_{FS}]$, where $\epsilon(\alpha,Y)=1$ and $\int_{X'}(\mu'^*\alpha-\{E\})^n=1/2$. 

Now we simply note that $$\int_{X'}\Omega^n=\int_{X'}(\mu'^*\alpha-c\{E\})^n\geq \int_{X'}(\mu'^*\alpha-\epsilon(\alpha,Y)\{E\})^n,$$ and so $$\int_{N_{Y|X}}\Omega^n\leq \int_X \omega^n-\int_{X'}(\mu'^*\alpha-\epsilon(\alpha,Y)\{E\})^n.$$ 

In fact this bound is easily seen to be sharp: for any $\epsilon>0$ there is a K\"ahler deformation $(\mathcal{X},\Omega)$ such that $$\int_{N_{Y|X}}\Omega^n\geq \int_X \omega^n-\int_{X'}(\mu'^*\alpha-\epsilon(\alpha,Y)\{E\})^n-\epsilon.$$

\subsection{Nonstandard deformations to the normal bundle}

If $\mathcal{X}'$ is the blow-up of $\mathcal{X}$ along some submanifold in $X'$, or more generally, if $\mathcal{X}'$ is a smooth modification of $\mathcal{X}$ with center contained in $X'$, then $\mathcal{X}'$ still contains $N_{Y|X}$ as part of the zero fiber. If we think of $\mathcal{X}_{core}:=\mathcal{X}\setminus X'$ as the core of the deformation of $X$ to the normal bundle of $Y$, then $\mathcal{X}$ and $\mathcal{X}'$ are just different compactifications of that core deformation. We will call such $\mathcal{X}'$ (nonstandard) deformations to the normal bundle of $Y$, and refer to $\mathcal{X}$ as the (standard) deformation to the normal bundle.

\begin{definition}
A pair $(\mathcal{X}',\Omega)$, where $\mathcal{X}'$ is a (possibly nonstandard) deformation of $X$ to the normal bundle of $Y$ and $\Omega$ is an $S^1$-invariant K\"ahler form on $\mathcal{X}$ such that $\Omega_{|X_1}=\omega$, will be called a \emph{K\"ahler deformation of $(X,\omega)$ to the normal bundle of $Y$}.
\end{definition}

The first main result of this paper says that with nonstandard K\"ahler deformations we can make sure that $\int_{N_{Y|X}}\Omega^n$ is arbitrarily close to $\int_X\omega^n$.

\begin{thmA} \label{ThmA}
For any $\epsilon>0$ there exists a K\"ahler deformation $(\mathcal{X}',\Omega)$ of $(X,\omega)$ to the normal bundle of $Y$ such that $$\int_{N_{Y|X}}\Omega^n\geq (1-\epsilon)\int_X\omega^n.$$  
\end{thmA}

Is it possible to deform $(X,\omega)$ so that all of its volume ends up in $N_{Y|X}$? Clearly, if we want $\Omega$ to be K\"ahler on some compact $\mathcal{X}'$ this is impossible, so instead we allow $\Omega$ to be a more general closed positive current. If $\Omega$ is an $S^1$-invariant closed positive $(1,1)$-current on $\mathcal{X}$ such that $\Omega_{|X_1}=\omega$ we call $(\mathcal{X},\Omega)$ a weak K\"ahler deformation of $(X,\omega)$ to the normal bundle of $Y$.

\begin{thmB} \label{ThmB}
There exists a weak K\"ahler deformation $(\mathcal{X},\Omega)$ of $(X,\omega)$ to the normal bundle of $Y$ such that $$\int_{N_{Y|X}}\Omega^n=\int_X\omega^n.$$
\end{thmB}
 
The weak K\"ahler deformation $(\mathcal{X},\Omega)$ in Theorem B is far from unique. However, we can construct a canonical weak K\"ahler deformation $(\mathcal{X}_{\overline{\mathbb{D}}},\Omega_{can})$ to the normal bundle of $Y$, where $\mathcal{X}_{\overline{\mathbb{D}}}:=\mu^{-1}(X\times \overline{\mathbb{D}})$. Apart from satisfying the volume equality $$\int_{N_{Y|X}}\Omega_{can}^n=\int_X\omega^n$$ it has have the additional property that $\Omega_{can}^{n+1}=0$. Hence it corresponds to a weak geodesic ray in the space of K\"ahler potentials with respect to $\omega$ equipped with the Mabuchi metric. The construction of $\Omega_{can}$ is described in Section \ref{Sec:can}.

\subsection{More general deformations}

Theorem A and B can be generalized in the following way.

Let $\mathcal{Z}$ be a smooth modification of $X\times \mathbb{P}^1$ with center contained in $X\times \{0\}$, and such that the zero-fiber is a SNC divisor, i.e. the union of smooth hypersurfaces $Z_0,...,Z_m$ with normal crossings. Let $D:=\cup_{i>0} Z_i$ and $Z:=Z_0\setminus D$.

\begin{definition}
A pair $(\mathcal{Z}',\Omega)$, where $\mathcal{Z}'$ is a smooth modification of $\mathcal{Z}$ with center contained in $D$, and $\Omega$ is an K\"ahler form on $\mathcal{Z}'$ such that $\Omega_{|X_1}=\omega$, will be called a K\"ahler deformation of $(X,\omega)$ to $Z$. If $\Omega$ is just closed and positive $(\mathcal{Z}',\Omega)$ is called a weak K\"ahler deformation of $(X,\omega)$ to $Z$. If the $\mathbb{C}^*$-action on $X\times \mathbb{P}^1$ lifts to $\mathcal{Z}$ we also demand that it lifts to $\mathcal{Z}'$ and that $\Omega$ is $S^1$-invariant with respect to this action.  
\end{definition}

\begin{thmAp}
For any $\epsilon>0$ there exists a K\"ahler deformation $(\mathcal{Z}',\Omega)$ of $(X,\omega)$ to $Z$ such that $$\int_{Z}\Omega^n\geq (1-\epsilon)\int_X\omega^n.$$ 
\end{thmAp}

\begin{thmBp}
There exists a weak K\"ahler deformation $(\mathcal{Z},\Omega)$ of $(X,\omega)$ to $Z$ such that $$\int_{Z}\Omega^n=\int_X\omega^n.$$ 
\end{thmBp}

If $Y$ is a singular subvariety of $X$, then by resolving singularities in the deformation to the normal cone of $Y$ we can get a smooth modification $\mathcal{Z}$ as above, together with an identification between a Zariski open and dense subset of $Z_0$ with a Zariski open and dense subset of the normal cone of $Y$. Hence Theorem A' shows that $(X,\omega)$ can be deformed so that almost all volume ends up on a regularization of the normal cone of $Y$.

\subsection{Big classes in $H^{1,1}(X,\mathbb{R})$}

The proofs of the above stated theorems rely on a general result about volumes of big cohomology classes, which is the other main result of this paper. First we recall some basic definitions regarding $(1,1)$ cohomology classes (note that we implicitly use the isomorphism between Dolbeault and Bott-Chern cohomology in the K\"ahler setting).

If $X$ is a compact K\"ahler manifold there are four convex cones in $H^{1,1}(X,\mathbb{R})$ corresponding to four different notions of positivity. First we have the \emph{K\"ahler cone} $\mathcal{K}$ consisting of the K\"ahler classes. This cone is open and its closure is called the \emph{nef cone} $\overline{\mathcal{K}}$. Then comes the \emph{pseudoeffective cone} $\mathcal{E}$, consisting of those classes that contain a closed positive current. This cone is closed, and its interior is called the \emph{big cone} $\mathcal{E}^{\circ}$. Thus a class is big iff it can be written as the sum of a pseudoeffective class and a K\"ahler class. 

When $X$ is projective these notions of positivity are consistent with those for holomorphic line bundles, i.e. a line bundle is ample/nef/big/pseudoeffective iff $c_1(L)$ is K\"ahler/nef/big/pseudoeffective.

Let $T$ be a closed positive current of bidegree $(1,1)$ on a compact K\"ahler manifold $(X,\omega)$, and let $\alpha$ be its cohomology class. Pick a smooth form $\theta$ in $\alpha$. By the $dd^c$-lemma $T$ can be written $dd^c u+\theta$, where $u$ is a quasi-psh function (i.e. locally the sum of a psh and a smooth function). A quasi-psh function $u$ such that $dd^c u+\theta$ is positive is called $\theta$-psh (the set of $\theta$-psh functions is denoted $PSH(X,\theta)$).

A closed current $T\in \alpha$ is called a \emph{K\"ahler current} if for some $\epsilon>0$, $T-\epsilon\omega$ is positive.

We say that $T=dd^cu+\theta$ has \emph{analytic singularities} if locally $u=a\ln(\sum_i|f_i|^2)+g$ where $f_i$ is a finite tuple of holomorphic functions and $g$ is smooth. Note that $T$ is smooth away from the analytic subset $E_T$ locally given by the common zero-set of $f_i$.

\begin{definition}
The intersection of the sets $E_T$ for all K\"ahler currents $T\in\alpha$ with analytic singularities is called the \emph{non-K\"ahler locus} of $\alpha$, denoted by $E_{nK}(\alpha)$, and is itself an analytic set. Its complement $K(\alpha)$ is called the \emph{K\"ahler locus}.
\end{definition}

When $X$ is projective and $\alpha=c_1(L)$ then $E_{nK}(\alpha)$ coincides with the augmented base locus of $L$.

A quantitative measure of bigness is provided by the volume function.

\begin{definition}
The volume $\textrm{vol}(\alpha)$ of a big class $\alpha$ is defined as the supremum of $\int_{X\setminus E_T} T^n$ for all K\"ahler currents $T\in\alpha$ with analytic singularities. If the class is not big we set the volume to be zero.
\end{definition}

If $\alpha$ is K\"ahler or nef $\textrm{vol}(\alpha)=\int_X\alpha^n,$ but this is not necessarily true when $\alpha$ is not nef.

If $\alpha$ and $\beta$ are big then it easy to see that $$\textrm{vol}(\alpha+\beta)\geq \textrm{vol}(\alpha)+\textrm{vol}(\beta).$$ Also, for any $c>0$: $\textrm{vol}(c\alpha)=c^n\textrm{vol}(\alpha)$. From these two properties it follows easily that the volume is continuous on the big cone. In fact, Boucksom proved in \cite{Bou} that it is continuous everywhere, i.e. that the volume goes to zero as one approaches the boundary of the big cone. 

There is also a notion of restricted volume along a submanifold $Y$.

\begin{definition}
The restricted volume $\textrm{vol}_{X|Y}(\alpha)$ of a big class $\alpha$ along a submanifold $Y$ of dimension $m$ which intersects $K(\alpha)$ is defined as the supremum of the integrals $\int_{Y\setminus E_T}T^m$ for all K\"ahler currents $T\in\alpha$ with analytic singularities. 
\end{definition} 
 
The same argument as for the ordinary volume shows it to be continuous on the open cone of big classes $\alpha$ such that $K(\alpha)$ intersects $Y$. The restricted volume is then extended to the whole $H^{1,1}(X,\mathbb{R})$ in the following way: if $K(\alpha+\epsilon[\omega])$ intersects $Y$ for all $\epsilon>0$, then $\textrm{vol}_{X|Y}(\alpha):=\lim_{\epsilon\to 0+}\textrm{vol}_{X|Y}(\alpha+\epsilon[\omega])$, otherwise $\textrm{vol}_{X|Y}(\alpha):=0$. 
 
When $X$ is projective and $\alpha=c_1(L)$ the volume of $\alpha$ coincides with the volume of $L$, and the same holds for the restricted volumes. 

We can now state our main result about big classes, which generalizes the corresponding result for the volume of line bundles proved independently by Boucksom-Jonsson-Favre \cite{BFJ} and Lazarsfeld-Musta\c{t}\u{a} \cite{LM}.

\begin{thmC} \label{Mainthm}
If $X$ is compact K\"ahler of dimension $n$, $\alpha\in H^{1,1}(X,\mathbb{R})$ is a big class, and $Y$ is a smooth hypersurface intersecting $K(\alpha)$, then $$\frac{d}{dt}_{|t=0}\textrm{vol}(\alpha+t\{Y\})=n\textrm{vol}_{X|Y}(\alpha).$$
\end{thmC}

A useful consequence is the following: 

\begin{corA} \label{Maincor}
If $Y_1,...,Y_m$ and $Z_1,...,Z_l$ are smooth hypersurfaces all of them intersecting $K(\alpha)$ and $$\sum_i\{Y_i\}=\sum_j\{Z_j\}$$ then $$\sum_i\textrm{vol}_{X|Y_i}(\alpha)=\sum_j\textrm{vol}_{X|Z_j}(\alpha).$$
\end{corA}

In particular, if $Y_{\tau}$ is a family of smooth hypersurfaces all intersecting $K(\alpha)$, then the restricted volume of $\alpha$ along $Y_{\tau}$ is independent of $\tau$. 

In fact we can generalize this a bit. The Lelong number $\nu_Y(\alpha)$ of $\alpha$ along $Y$ is defined as the infimum of the Lelong numbers $\nu_Y(T)$ over all closed positive currents $T\in \alpha$ (see Section \ref{Sec:lelong}). It is not hard to show that if $\nu_Y(\alpha)>0$ then both $\frac{d}{dt}_{|t=0}\textrm{vol}(\alpha+t\{Y\})$ and $\textrm{vol}_{X|Y}(\alpha)$ are zero. Thus in Corollary C we can also allow hypersurfaces with positive Lelong number.

\subsection{On the proof of Theorem C}

We will use the deformation $\mathcal{X}$ of $X$ to the normal bundle of $Y$, already discussed above, equipped with the classes $$\beta=(\pi_X\circ\mu)^*\alpha+b\{X_0\}-c\{\mathcal{E}\},$$ with $b>c>0$. We choose $c$ so that $\alpha-c\{Y\}$ is big.

The proof is divided into four steps.

\subsection{Step 1: Volume formulas}

The first thing we want to do is to express the volume and restricted volumes of $\beta$ in terms of the volumes and restricted volumes of $\alpha-t\{Y\}$:

\begin{theorem} \label{volcalc2}
\begin{enumerate}
\item $$\textrm{vol}(\beta)=(n+1)\left((b-c)\textrm{vol}(\alpha)+\int_0^c\textrm{vol}(\alpha-t\{Y\})dt\right),$$
\item $$\textrm{vol}_{\mathcal{X}|X_{\tau}}(\beta)=\textrm{vol}(\alpha), \qquad{} \forall \tau\neq 0,$$
\item $$\textrm{vol}_{\mathcal{X}|X'}(\beta)=\textrm{vol}(\alpha-c\{Y\}),$$
\item $$\textrm{vol}_{\mathcal{X}|\mathcal{E}}(\beta)=n\int_0^c\textrm{vol}_{X|Y}(\alpha-t\{Y\})dt.$$
\end{enumerate}
\end{theorem}

To prove these formulas we will need a couple of things. First, recall that a current $T=dd^cu+\theta$ is said to have \emph{minimal singularities} if $u\geq v-C$ for any $v\in PSH(X,\theta)$. A key fact proved in \cite{BEGZ} is that if $T$ has minimal singularities then $$\int_X T^n=\textrm{vol}(\alpha),$$ where $T^n$ stands for the non-pluripolar Monge-Amp\`ere measure of $T$ (see Section \ref{Sec:MA}). 

A common way to find currents with minimal singularities is via envelopes. If $g$ is a smooth function on $X$ then the associated envelope $u$ is defined as $$u:=\sup\{v\leq g: v\in PSH(X,\theta)\}.$$ The current $T=dd^cu+\theta$ will be positive and have minimal singularities.

These envelopes have an additional property which makes them extremely useful in calculating volumes. Namely, if $D:=\{u=g\}$ is the contact set, then $$T^n=\mathbbm{1}_D(dd^cg+\theta)^n.$$ When $[\theta]$ is K\"ahler this is due to Berman \cite{Ber13}, while the big case was proved by Di Nezza-Trapani \cite{DNT}. 

The second important tool we use is the partial Legendre transform due to Kiselman \cite{Kis}. Recall that the Legendre transform is a transform on the space of convex functions. If $u(\tau)$ is subharmonic in a complex variable $\tau$ and only depends on $|\tau|$, then $u$ is convex in the variable $y:=\ln|z|^2$. When $u(z,\tau)$ is $\pi_X^*\theta$-psh on $X\times \mathbb{C}^*$ and $u(z,\tau)=u(z,|\tau|)$, performing a Legendre transform along each orbit we get the partial Legendre transform, and the effect is a decomposition of $u$ into a concave family of $\theta$-psh functions $\hat{u}_{\lambda}$ on $X$.

In Section \ref{Sec:trans} we prove a general theorem, Theorem \ref{Thmlocal}, possibly of independent interest, which expresses the Monge-Amp\`ere measure of $u$ in terms of the Monge-Amp\`ere measures of the transform $\hat{u}_{\lambda}$: $$(\pi_X)_*(dd^cu+\pi_X^*\theta)^{n+1}=(n+1)\int_{\lambda}(dd^c\hat{u}_{\lambda}+\theta)^n d\lambda.$$

To prove Theorem \ref{volcalc2} we will use an envelope to construct a current $T\in \beta$ with minimal singularities. We will then use the partial Legendre transform to decompose it into currents on $X$ and analyse these. Finally Theorem \ref{Thmlocal} will allow us express the volume and restricted volume of $\beta$ in terms of the volume and restricted volume of $\alpha-t\{E\}$.

\subsection{Step 2: The key equality}

We then want to prove the following key equality:

\begin{theorem} \label{keyeq}
\begin{eqnarray} \label{aa}
\textrm{vol}_{\mathcal{X}|X'}(\beta)+\textrm{vol}_{\mathcal{X}|\mathcal{E}}(\beta)=\textrm{vol}(\alpha).
\end{eqnarray}
\end{theorem}

Combined with Theorem \ref{volcalc2} the equality will give us that $$\textrm{vol}(\alpha)=\textrm{vol}(\alpha-c\{Y\})+n\int_0^c\textrm{vol}_{X|Y}(\alpha-t\{Y\})dt.$$ As $\textrm{vol}_{X|Y}(\alpha-t\{Y\})$ is continuous in $t$ at $t=0$, this will imply that $$\frac{d}{dt}_{|t=0}\textrm{vol}(\alpha+t\{Y\})=n\textrm{vol}_{X|Y}(\alpha),$$ i.e. Theorem C.  

One inequality of (\ref{aa}) is straight-forward; the difficulty lies in establishing $$\textrm{vol}_{\mathcal{X}|X'}(\beta)+\textrm{vol}_{\mathcal{X}|\mathcal{E}}(\beta)\geq \textrm{vol}(\alpha).$$

We thus come to the third step of the proof.

\subsection{Step 3: Measure control}

We let $\tilde{\theta}$ be a smooth form in $\beta$ and define $$u:=\sup\{v\leq 0: v\in PSH(\mathcal{X},\tilde{\theta})\}.$$ Then $T:=dd^c u+\tilde{\theta}$ has minimal singularities and $\int_{X_{\tau}}T^n= \textrm{vol}(\alpha)$ for all $|\tau|>0$. It can be proved that the measures $T_{|X_{\tau}}^n$ converge weakly to $T_{|X'}^n$ and $T_{|\mathcal{E}}^n$ locally near $(X'\setminus E)\cap K(\beta)$ and $(\mathcal{E}\setminus E)\cap K(\beta)$, and this implies the easy inequality. But to prove the other inequality we have to make sure that none of the mass of $T_{|X_{\tau}}^n$ dissappears into the set of singularities $E\cup (X_0\cap E_{nK}(\beta))$, as this would lead to $$\textrm{vol}_{\mathcal{X}|X'}(\beta)+\textrm{vol}_{\mathcal{X}|\mathcal{E}}(\beta)< \textrm{vol}(\alpha).$$ 

Since $T$ was defined via an envelope $u$ we have that 
\begin{eqnarray*} \label{est1}
T^{n+1}=\mathbbm{1}_{D}\tilde{\theta}^{n+1},
\end{eqnarray*}
where $D:=\{u=0\}$ is the contact set, while for $T^n_{|X_{\tau}}$ we have the inequality
\begin{eqnarray} \label{est2}
T^n_{|X_{\tau}}\geq \mathbbm{1}_{X_{\tau}\cap D}\tilde{\theta}^n.
\end{eqnarray}
Suppose that for some sequence $\tau_k\to 0$ we knew that 
\begin{equation} \label{eq100}
\lim_{k\to \infty}\int_{X_{\tau_k}\cap D}\tilde{\theta}^n=\textrm{vol}(\alpha).
\end{equation} 
For any $\epsilon>0$ we can find an open neighbourhood $U$ of the set of singularities $E\cup (X_0\cap E_{nK}(\beta))$ such that for all $|\tau|>0$: 
\begin{equation} \label{eq101}
\int_{X_{\tau}\cap U}|\tilde{\theta}^n|<\epsilon.
\end{equation}
The local convergence of $T_{|X_{\tau}}^n$ away from $U$ together with the estimates (\ref{est2}), (\ref{eq100}) and (\ref{eq101}) then implies that $$\int_{X'}T^n+\int_{\mathcal{E}}T^n\geq \liminf_{k\to\infty}\int_{(X_{\tau_k}\cap D)\setminus U}\tilde{\theta}^n\geq \textrm{vol}(\alpha)-\epsilon,$$ and since $\epsilon$ was arbitrary we would be done. Unfortunately (\ref{eq100}) seems unlikely to hold.  

Instead our strategy will be to construct a sequence of envelope currents $T_k$ with contact sets $D_k$ together with a sequence $\tau_k\to 0$ such that 
\begin{equation} \label{est4}
\lim_{k\to \infty}\int_{X_{\tau_k}\cap D_k}\tilde{\theta}^n=\textrm{vol}(\alpha).
\end{equation}
We then also need to show that the measures ${T_k^n}_{|X_{\tau_k}}$ converge locally to some measures  $\mu_{X'}$ and $\mu_{\mathcal{E}}$ away from the set of singularities, and that $$\int_{X'}\mu_{X'}\leq \textrm{vol}_{\mathcal{X}|X'}(\beta)$$ and $$\int_{\mathcal{E}}\mu_{\mathcal{E}}\leq \textrm{vol}_{\mathcal{X}|\mathcal{E}}(\beta).$$ If we can do all that, then by the same arguments as above, this will establish the desired inequality, and hence conclude the proof of Theorem C. 

To find currents $T_k$ with the desired properties we will elaborate on a technique introduced in \cite{WN2}. The trick is to consider a sequence of smooth functions $g_k$ that converge to the singular function $\ln|\tau|^2-\ln(1+|\tau|^2)$. We let $u_k$ be the corresponding envelopes, $T_k:=dd^cu_k+\tilde{\theta}$ and $D_k:=\{u_k=g_k\}$.

To establish (\ref{est4}) we will argue in the following way. By construction $T_k\to T_{\infty}+[X_0]$ where $T_{\infty}$ is a closed positive current in $\beta-\{X_0\}$. Using ideas from \cite{WN2} together with Theorem \ref{volcalc2} we prove that for any $\epsilon>0:$ $$\liminf_{k\to\infty}\int_{|\tau|\leq \epsilon}T_k^{n+1}\geq (n+1)\textrm{vol}(\alpha).$$ We also have that $$T_k^{n+1}=\mathbbm{1}_{D_k}(\tilde{\theta}+dd^cg_k)^{n+1}\leq \mathbbm{1}_{D_k}(\tilde{\theta}^n+(n+1)\tilde{\theta}^n\wedge \omega_k),$$ where $\omega_k:=\omega_{FS}+dd^cg_k.$ Together this implies that $$\liminf_{k\to \infty}\int_{\mathbb{P}^1}\left(\int_{X_{\tau}\cap D_k}\tilde{\theta}^n\right)\omega_k\geq\textrm{vol}(\alpha).$$ Finally, since $\omega_k$ are probability measures on $\mathbb{P}^1$ that converge to $\delta_0$, and for each $\tau$ $$\int_{X_{\tau}\cap D_k}\tilde{\theta}^n\leq \textrm{vol}(\alpha),$$ this means that one can find a sequence $\tau_k$ as in (\ref{est4}).

\subsection{Step 4: Convergence}

We then come to the question of convergence of the measures ${T_k^n}_{|X_{\tau_k}}$. 

The basic idea is as follows. We will choose $g_k$ so that $g_k(\tau)\approx g(e^{k/2}\tau)-k$ which will give us a Laplacian bound on $g_k$ of order $e^k$. Combined with a key regularity result for envelopes due to Berman \cite{Ber13} this will yield a Laplacian bound on the envelopes $u_k$ of order $e^k$. The Laplacian bound together with a supremum bound then will result in a local Lipschitz bound, which is enough to establish the local convergence of ${T_k^n}_{|X_{\tau_k}}$ as long as $e^{k/2}\tau_k\to 0$. Luckily, Step 3 will allow us to choose $\tau_k$ so that this holds.

There is a technical issue though, namely that Berman's regularity result demands that the cohomology class one is working with is K\"ahler, and not only big as in our case. This will be handled by using approximate Zariski decompositions of $\beta$ on suitible modifications of $\mathcal{X}$, and then passing to the limit (see Section \ref{Sec:conv}). 

\subsection{On the proofs of Theorem A and B}

Let $\alpha:=[\omega]$. For notational simplicity assume that $Y$ is a hypersurface. Recall that the supremum of $t$ such that $\alpha-t\{Y\}$ is K\"ahler was called the Seshadri constant $\epsilon(\alpha,Y)$. The related \emph{pseudoeffective threshold} $\delta(\alpha,Y)$ is defined as the supremum of $t$ such that is $\alpha-t\{Y\}$ is pseudoeffective. By continuity, the volume of $\alpha-t\{Y\}$ tends to zero as $t\to \delta(\alpha,Y)$. Thus, given $\epsilon>0$ we can pick a $c$ such that $$0<\textrm{vol}(\alpha-c\{Y\})<\epsilon\textrm{vol}(\alpha).$$ Also pick some $b>c$ and let $$\beta=(\pi_X\circ\mu)^*\alpha+b\{X_0\}-c\{\mathcal{E}\}.$$

Combining Theorem \ref{volcalc2} and Corollary C we get that $$\textrm{vol}_{\mathcal{X}|\mathcal{E}}(\beta)=\textrm{vol}(\alpha)-\textrm{vol}(\alpha-c\{Y\})\geq (1-\epsilon)\textrm{vol}(\alpha).$$ By the definition of restricted volume this means that one can find a K\"ahler current $T\in \beta$ with analytic singularities such that $$\int_{N_{Y|X}\setminus E_T}T^n>(1-2\epsilon)\textrm{vol}(\alpha).$$ In fact $E_{nK}(\beta)\subset X'$, and wlog we can assume that $E_T=E_{nK}(\beta)$. There is then a smooth modification $\pi':\mathcal{X}'\to \mathcal{X}$ with center contained in $X'$ such that $\pi'^*T$ has divisorial singularities, i.e. $\pi'^*T=\Omega+\sum a_i[E_i]$. Thus $\Omega'$ is a semipositive form on $\mathcal{X}'$ such that $[\Omega'_{|X_1}]=\alpha$ and $\int_{N_{Y|X}}\Omega'>(1-2\epsilon)\textrm{vol}(\alpha)$. Given this it is easy to find an $S^1$-invariant K\"ahler form $\Omega$ such that $\Omega_{|X_1}=\omega$ and $$\int_{N_{Y|X}}\Omega>(1-2\epsilon)\textrm{vol}(\alpha).$$

To prove Theorem B we let $c=\delta(\alpha,Y)$ and pick an $\Omega$ with minimal singularities. It follows from the above calculations that $$\int_{N_{Y|X}}\Omega^n=\textrm{vol}(\alpha),$$ and it is again not hard to see that $\Omega$ can be chosen to be $S^1$-invariant and such that $\Omega_{|X_1}=\omega$.

The proof of Theorem A' and B' are similar.

\subsection{Related work} 

\subsubsection{Differentiability of volume}

When $X$ is projective the \emph{Neron-Severi space} $NS(X,\mathbb{R})\subseteq H^{1,1}(X,\mathbb{R})$ is the subspace generated by Chern classes of holomorphic line bundles. As was already mentioned above Boucksom-Jonsson-Favre \cite{BFJ} and Lazarsfeld-Musta\c{t}\u{a} \cite{LM} independently proved the line bundle version of Theorem C. Since any class in $NS(X,\mathbb{R})$ can be written as the difference of two ample classes, and any ample class can be approximated by rational multiples of smooth divisors, this implied that the volume restricted to the big cone in $NS(X,\mathbb{R})$ is differentiable, indeed $C^1$. Boucksom-Jonsson-Favre \cite{BFJ} also proved that $$\frac{d}{dt}_{|t=0}\textrm{vol}(\alpha+t\gamma)=n\langle \alpha^{n-1} \rangle \cdot \gamma,$$ where $\langle \alpha^{n-1} \rangle$ is a cohomology class called the positive selfintersection of $\alpha$. For more on the volume of line bundles see \cite{Laz04}. 

In \cite{BDPP} Boucksom-Demailly-P\u aun-Peternell conjectured that for any compact K\"ahler manifold the volume is differentiable on the full big cone $\mathcal{E}^{\circ}\subseteq H^{1,1}(X,\mathbb{R})$, and that the derivative is given by $$\frac{d}{dt}_{|t=0}\textrm{vol}(\alpha+t\gamma)=n\langle \alpha^{n-1} \rangle \cdot \gamma.$$ The special case of this conjecture when $X$ is projective was proved by the author in \cite{WN2}.

In \cite{ELMNP} Ein-Lazarsfeld-Musta\c{t}\u{a}-Nakamaye-Popa proved that the union of subvarieties along which a big line bundle has zero restricted volume is equal to its augmented base locus. In \cite{CT} Collins-Tosatti proved the same statement for nef classes on compact K\"ahler manifolds, but the big case is still open.

\subsubsection{Test configurations and geodesic rays}

The deformation spaces $\mathcal{X}$ appearing in this paper are examples of so-called \emph{test configurations}, which play a central role in the famous Yau-Tian-Donaldson conjecture (see e.g. \cite{Sze14} and references therein).

A test configuration of a smooth polarized projective variety $(X,L)$ is a normal polarized projective variety $(\mathcal{X},\mathcal{L})$ together with a $\mathbb{C}^*$-action on $\mathcal{X}$ lifting to $\mathcal{L}$ and a flat $\mathbb{C}^*$-equivarient map $\pi: \mathcal{X}\to \mathbb{P}^1$ such that $\pi^{-1}(\mathbb{P}^1\setminus \{0\},\mathcal{L})$ is $\mathbb{C}^*$-equivariantly isomorphic to $(X\times \mathbb{P}^1\setminus \{0\},\pi_X^*L^r)$ for some number $r$.

Let $\omega$ be a K\"ahler form in $c_1(L)$. Phong-Sturm showed in \cite{PS} that to any test configuration one can associate a weak geodesic ray in the space of K\"ahler potentials with respect to $\omega$ (or rather its completion). A weak geodesic ray can be interpreted as an $S^1$-invariant closed positive current $\Omega$ on $X\times \mathbb{D}^*$ such that $\Omega_{|X_1}=\omega$ and $\Omega^{n+1}=0$. In fact $\Omega$ will extend over the central fiber of the test configuration, and will represent the first Chern class of $\mathcal{L}$.

It also makes sense to consider test configurations where $\mathcal{L}$ is not ample but big. It follows from the work in \cite{RWN1} that with $\mathcal{L}$ big we still get a weak geodesic ray in the sense above. When $X$ is projective and $\omega\in c_1(L)$ the canonical deformation $\Omega_{can}$ of this paper corresponds (up to a trivial change) to the ray associated to $(\mathcal{X},\mathcal{L})$ where $\mathcal{X}$ is the deformation to the normal cone of $Y$ and $\mathcal{L}:=\pi_X^*L\otimes \pi_{\mathbb{P}^1}(L_H)^{N}\otimes L_{\mathcal{E}}^{-M}$, where $L_H$ is the hyperplane line bundle on $\mathbb{P}^1$, $L_{\mathcal{E}}$ is the line bundle on $\mathcal{X}$ corresponding to the divisor $\mathcal{E}$, and $N,M$ are natural numbers so that $N>M>\delta(\alpha,Y)$.

If instead of a smooth polarized variety we have a compact K\"ahler manifold $X$ with a K\"ahler class $\alpha$ one can mimick the definitions above, replacing the ample line bundle $\mathcal{L}$ with a K\"ahler class $\beta$, as is done in \cite{DR,SD}. Donaldson-Futaki invariants can then be defined as intersection numbers \cite{DR,SD}, giving risee to a notion of K-stability and a formulation of the Yau-Tian-Donaldson conjecture in the transcendental setting. It was proved independently by Dervan-Ross \cite{DR} and Sj\"ostr\"om Dyrefelt \cite{SD} that the existence of a cscK metric implies K-semistability.

As in the projective case test configurations give rise to weak geodesic rays, even when $\beta$ is not K\"ahler but big, as is often the case in this paper.

\subsubsection{The non-archimedean Monge-Amp\`ere equation}

Theorem A' is related to the work of Boucksom-Jonsson-Favre \cite{BFJ15, BJ} on the non-archimedean Monge-Amp\`ere equation, and hence also the variational approach to the Yau-Tian-Donaldson conjecture \cite{BBJ}. 

One motivation behind the work of Boucksom-Jonsson-Favre is to complete the space of test configurations.

Given a polarized projective variety $(X,L)$ over say the trivially valued field $\mathbb{C}$ there is an associated polarized Berkovich space $(X^{an},L^{an})$, called the analytification of $(X,L)$. A test configuration $(\mathcal{X}',\mathcal{L})$ of $(X,L)$ induces a continuous function on $X^{an}$ and is thought of as a positive metric on $L^{an}$. The components $D_i$ of the zero-fiber of $\mathcal{X}'$ correspond to points $x_i$ in $X^{an}$, and the non-archimedean Monge-Amp\`ere measure of the metric associated to $(\mathcal{X}',\mathcal{L})$ is defined as the atomic probability measure $$\frac{1}{(L^n)}\sum_i(\mathcal{L}^n \cdot D_i)\delta_{x_i}.$$ The space of singular positive metrics on $L^{an}$ is then defined as the set of decreasing limits of positive metrics. The Monge-Amp\`ere operator can be extended to the space of continuous singular metrics \cite{CLD} and even to the larger space of finite energy singular metrics \cite{BJ}. Using the results of \cite{BDPP} on the orthogonality of Zariski decomposition Boucksom-Favre-Jonsson proved a Calabi-Yau theorem\cite{BFJ15,BJ} in this setting, i.e. the Monge-Amp\`ere operator is a bijection between the space of finite energy singular positive metrics and the set of finity energy Radon probability measures on $X^{an}$. 

If $X$ is projective and $\omega\in c_1(L)$  for some ample line bundle, then the cohomology class $\beta$ on $\mathcal{X}$ in the proof of Theorem A is similarly the first Chern class of a big line bundle $\mathcal{L}$ on $\mathcal{X}$, which corresponds to a singular positive metric on $L^{an}$. The proof of Theorem A shows that the non-archimedean Monge-Amp\`ere measure of this metric is equal to the Dirac measure at the point in $X^{an}$ corresponding to $\mathcal{E}$. In Theorem A' we instead have a Dirac measure at the point corresponding to $Z_0$.   

\subsubsection{The partial Legendre transform}

The partial Legendre transform of plurisubharmonic functions with symmetry was introduced by Kiselman in \cite{Kis}, where the crucial minimum principle is proved. It has later been used to great effect e.g. by Demailly in \cite{Dem92} to find regularizations of plurisubharmonic functions. In \cite{RWN1} the partial Legendre transform was used, in a similar way as here, to analyse quasi-psh functions on test configurations, in particular geodesics. That technique has later been used by e.g. Darvas-Rubinstein \cite{DR} to study geodesic segments, and Darvas-Xia \cite{DX} to study limits of test configurations.

The partial Legendre transform was also used in \cite{WN3} to prove monotonicity of Monge-Amp\`ere masses (Theorem \ref{mono}). Key there was an expression of the Monge-Amp\`ere in terms of the Monge-Amp\`ere of the transforms, which here is generalized in Theorem \ref{Thmlocal}. This will be further explored in the forthcoming paper \cite{BWN}.

\subsubsection{Canonical tubular neighbourhoods}

In the paper \cite{RWN2} Ross and the author showed how to construct a canonical tubular neighbourhood of any submanifold $Y$ of a K\"ahler manifold $(X,\omega)$. This relied on finding a canonical $S^1$-invariant closed positive form $\Omega$ in a neighbourhood $U\subseteq \mathcal{X}_{\overline{\mathbb{D}}}$ of $\mu^{-1}(Y\times \mathbb{D})$. $\Omega$ was found by locally solving the homogeneous Monge-Amp\`ere equation with boundary data $\Omega_{|U_1}=\omega_{|U_1}$, and then showing that the local solutions agreed on overlaps. $\Omega$ then gave rise to a foliation of $U$, and the corresponding map from $U_1$ to $N_{Y|X}\cap U$ gave the canonical tubular neighbourhood of $Y$, specially adapted to the K\"ahler form $\omega$. This is very much connected to the canonical weak K\"ahler deformation $(\mathcal{X}_{\overline{\mathbb{D}}},\Omega_{can})$ of this paper, since ${\Omega_{can}}_{|U}$ coincides with the canonical form $\Omega$ from \cite{RWN2}. 

\subsubsection{The Hele-Shaw flow}

The canonical deformation $(\mathcal{X}_{\overline{\mathbb{D}}},\Omega_{can})$ in the special case where $X=\mathbb{P}^1$ and $Y=\{0\}$ was considered already in \cite{RWN3} (see also \cite{RWN18}). There it was shown that via the partial Legendre transform the deformation was equivalent to the Hele-Shaw flow, with $\omega$ encoding the permeability of the medium.

\subsubsection{Canonical growth conditions}

The special case of Theorem \ref{ThmB} when $X$ is projective, $\omega$ is the curvature form of a positive metric on an ample line bundle, and $Y$ is a point, was used in \cite{WN1} to analyse the so-called canonical growth condition.  

\begin{ackn}
I want to thank Bo Berndtsson, Robert Berman and Julius Ross for many enlightening discussions on related topics over the years. 
\end{ackn}

\section{Preliminaries} \label{Sec:prel}

Let $(X,\omega)$ be a compact K\"ahler manifold of complex dimension $n$, $\alpha\in H^{1,1}(X,\mathbb{R})$ a big class and $\theta$ a smooth form in $\alpha$.

\subsection{Lelong numbers} \label{Sec:lelong}

If $z_i$ are local holomorphic coordinates centered at a point $p\in X$, then the Lelong number $\nu_p(u)$ of $u\in PSH(X,\theta)$ at $p$ is defined as the supremum of $\lambda$ such that locally $u\leq \lambda\ln||z||^2+C$. If $Y$ is a subvariety then the Lelong number $\nu_Y(u)$ of $T$ along $Y$ is defined as the infimum $$\nu_Y(u):=\inf\{\nu_p(u): p\in Y\}.$$  If $Y$ is a divisor then $\nu_Y(u)\geq \lambda$ iff $dd^cu+\theta-\lambda[Y]$ is positive.

If $T=dd^cu+\theta$ we let $\nu_p(T):=\nu_p(u)$ and $\nu_Y(T):=\nu_Y(u)$.  

\subsection{Regularization of quasi-psh functions} 

A fundamental result we need to mention is Demailly's regularization theorem \cite[Thm. 1.1]{Dem92}:

\begin{theorem}
If $u\in PSH(X,\theta)$ then there is a sequence of functions $u_j\in PSH(X,\theta+\omega/j)$ with analytic singularities decreasing to $u$ and such that the Lelong numbers of $u_j$ increases to those of $u$. 
\end{theorem}

The full statement of the result also includes more precise control of the differences between $u_j-u$, but this will not be needed here.

As a consequence any big class contains a K\"ahler current with analytic singularities.

\subsection{Non-pluripolar Monge-Amp\`ere measures} \label{Sec:MA}

If $u\in PSH(X,\theta)$ is locally bounded on an open set $U$ then by the work of Bedford-Taylor \cite{BT} $(dd^cu+\theta)^m$ is a well-defined closed positive current on $U$. When $m=n$ we thus get a positive measure on $U$, called the Monge-Amp\`ere measure of $U$ (with respect to $\theta$), also denoted $MA_{\theta}(u)$.

The Monge-Amp\`ere measure does not charge pluripolar sets. In \cite{BT} Bedford-Taylor also established the following absolutely crucial continuity property for the Monge-Amp\`ere operator:

\begin{theorem} \label{BT}
If $u_j$ is a sequence of locally bounded $\theta$-psh functions decreasing or increasing a.e. to a locally bounded $\theta$-psh function $u$, then $MA_{\theta}(u_j)$ converge weakly to $MA_{\theta}(u)$.
\end{theorem} 

Another important property of the Monge-Amp\`ere measure established by Bedford-Taylor \cite{BT2} is that it is local in the plurifine topology. The plurifine topology is defined as the coarsest topology making quasi-psh functions continuous. That the Monge-Amp\`ere operator is local with respect to this topology means in particular that if $u,v\in PSH(X,\theta)$ and $u=v$ on some plurifine open set $O$ then $$\mathbbm{1}_O MA_{\theta}(u)=\mathbbm{1}_O MA_{\theta}(v).$$

Let $u\in PSH(X,\theta)$ be unbounded, and let $U$ be an open set on which there is a bounded $\theta$-psh function $v$. Then for all $j$, $u_j:=\max(u,v-j)$ is $\theta$-psh and bounded, thus $MA_{\theta}(u_j)$ is a well-defined measure on $U$. Note that by locality $\mathbbm{1}_{\{u>v-j\}}MA_{\theta}(u_j)$ is increasing, and the non-pluripolar Monge-Amp\`ere measure of $u$ is defined on $U$ as the limit of this as $j\to \infty$. That it does not depends on the particular choice of $v$ also follows from locality, and so it defines a measure on $X.$ A priori one could have ended up with something not locally finite, but it was shown in \cite{BEGZ} that, thanks to $X$ being compact K\"ahler, it is a finite measure on $X$. It is called the non-pluripolar Monge-Amp\`ere measure of $u$ (with respect to $\theta$), and also denoted $MA_{\theta}(u)$.   

The non-pluripolar Monge-Amp\`ere measure does not charge pluripolar sets, and is local in the plurifine topology. It is however not continuous under decreasing sequences.

On $PSH(X,\theta)$ there is a natural partial order, namely we write $u\succeq v$ if $u\geq v-C$ for some constant $C$. We then say that $u$ is less singular than $v$. A useful fact, proved in \cite{WN3}, is that the total Monge-Amp\`ere mass is monotone with respect to this partial order:

\begin{theorem} \label{mono}
If $u,v\in PSH(X,\theta)$ and $u\succeq v$ then $$\int_X MA_{\theta}(u)\geq \int_X MA_{\theta}(v).$$
\end{theorem}
The most important case when $u$ and $v$ are locally bounded outside a closed pluripolar set was already proved in \cite{BEGZ}.

If $u\in PSH(X,\theta)$ is maximal with respect to this partial order we say that $u$ has minimal singularities. It follows from Theorem \ref{mono} that if $u$ has minimal singularities then $$\int_X MA_{\theta}(u)\geq \textrm{vol}(\alpha).$$ On the other hand it follows from Demailly's regularization theorem that for any $\epsilon>0$ we can find an $\theta+\epsilon\omega$-psh function $u_{\epsilon}\geq u$ such that $dd^c u_{\epsilon}+\theta+\epsilon\omega$ is a K\"ahler current with analytic singularities. We thus get that 
\begin{eqnarray*}
\int_X MA_{\theta}(u)\leq \int_X MA_{\theta+\epsilon\omega}(u)\leq \int_X MA_{\theta+\epsilon\omega}(u_{\epsilon})\leq \textrm{vol}(\alpha+\epsilon[\omega]).
\end{eqnarray*} 

Since the volume function is continuous we get the result, already established in \cite{BEGZ}, that $$\int_X MA_{\theta}(u)=\textrm{vol}(\alpha)$$ whenever $u\in PSH(X,\theta)$ has minimal singularities. 

Exactly the same argument shows that if $u\in PSH(X,\theta)$ has minimal singularities and $K(\alpha)$ intersects $Y$ we have that $$\int_Y MA_{\theta_{|Y}}(u_{|Y})=\textrm{vol}_{X|Y}(\alpha).$$

\subsection{Envelopes}

Recall from the introduction that if $\theta\in \alpha$ is smooth and $g$ is a continuous function on $X$ the corresponding envelope $u$ is defined as $$u:=\sup\{v\leq g: v\in PSH(X,\theta)\}.$$

Let us record the following general result on Monge-Amp\`ere measures on contact sets due to Di Nezza-Trapani \cite{DNT}: 

\begin{theorem} \label{DNT}
If $u\in PSH(X,\theta)$ and $u\leq g$ where $g$ is a continuous function with bounded distributional Laplacian, then $$\mathbbm{1}_{\{u=g\}}MA_{\theta}(u)=\mathbbm{1}_{\{u=g\}}(dd^cg+\theta)^n.$$
\end{theorem}

If $u$ is the envelope with respect to $g$ then it is a classical fact due to Bedford-Taylor \cite{BT} that $\mathbbm{1}_{\{u<g\}}MA_{\theta}(u)=0$, and so combined with Theorem \ref{DNT} we get $$MA_{\theta}(u)=\mathbbm{1}_{\{u=g\}}(dd^cg+\theta)^n.$$ Let us also note that if $Y$ is a submanifold of dimension $m$ then Theorem \ref{DNT} applied to $u_{|Y}$ and $g_{|Y}$ gives that $$\mathbbm{1}_{\{u=g\}\cap Y}MA_{\theta_{|Y}}(u_{|Y})=\mathbbm{1}_{\{u=g\}\cap Y}(dd^cg+\theta)^m$$ and hence $$MA_{\theta_{|Y}}(u_{|Y})\geq \mathbbm{1}_{\{u=g\}\cap Y}(dd^cg+\theta)^m.$$

We will also need the following regularity result for envelopes due to Berman \cite{Ber13}. Note that here the reference form $\theta$ is supposed to be K\"ahler.

\begin{theorem} \label{ThmBerman}
Let $(X,\theta)$ be compact K\"ahler, $g$ a smooth function and $\psi$ a quasi-psh function with analytic singularities, and let $u$ be the corresponding envelope $$u:=\sup\{v\leq g-\psi: v\in PSH(X,\theta)\}.$$ Then we have the Laplacian estimate $$0\leq n+\Delta_{\theta}u\leq ((C+1)n+\sup_X(\Delta_{\theta}g))e^{B(g-\psi-\inf_X(g-\psi))},$$ where $C$ is a constant such that $dd^c\psi\geq -C\theta$ and $-B$ is a negative lower bound of the holomorphic bisectional curvature of $\theta$. 
\end{theorem}

\begin{proof}
In \cite{Ber13} this is done without the quasi-psh function $\psi$, so let us just explain how to reduce it to the smooth case proved by Berman.

Since $u$ is bounded and $g$ is smooth $u-g\leq R$ for some constant $R$. It follows that $u-g\leq \min(-\psi,R)=-\max(\psi,-R)$ and so $u\leq g-\max(\psi,-R)$. Let $\max_{reg}$ be a regularized $max$-function such that $\max_{reg}(x,y)=\max(x,y)$ when $|x-y|\geq 1$ say, and let $\tilde{\psi}:=\max_{reg}(\psi,-R-1)$. Then $\tilde{\psi}$ is smooth, $dd^c\tilde{\psi}\geq -C\theta$ and by definition $\tilde{\psi}\leq \max(\psi,-R)$. We thus get that $u\leq g-\tilde{\psi}$, and hence that $$u\leq \sup\{v\leq g-\tilde{\psi}:v\in PSH(X,\theta)\}.$$ On the other hand, since $\tilde{\psi}\geq \psi$ we have that $$\sup\{v\leq g-\tilde{\psi}: v\in PSH(X,\theta)\}\leq u$$ and hence we get equality. This shows that without loss of generality $\psi$ can be assumed to be smooth.

\end{proof}

\section{The partial Legendre transform} \label{Sec:trans}

We start by recalling the Legendre transform of convex functions on $\mathbb{R}$. Let $f: \mathbb{R}\to \mathbb{R}\cup\{\infty\}$ be a convex function. The Legendre transform $\hat{f}$ is a new convex function on $\mathbb{R}$ defined by $$\hat{f}(\lambda):=\sup_{y\in \mathbb{R}}\{\lambda y-f(y)\}.$$ Note that if $\lambda$ is a subgradient of $f$ at the point $y_{\lambda}$ then $\hat{f}(\lambda)=\lambda y_{\lambda}-f(y_{\lambda})$. In particular, if $f$ is smooth and strictly convex then $y_{\lambda}$ depends smoothly on $\lambda$ and thus $\hat{f}$ is also smooth on the interval where it is finite.

One also sees that if we take the Legendre transform of $\hat{f}$ we get the supremum of all affine functions bounded from above by $f$. Since $f$ is convex this is precisely $f$, except possibly on the boundary of the interval where $f$ is finite. On the left limit point $a$ (if it exists) we then get that $\hat{\hat{f}}(a)=\lim_{y\to a+}f(y)$ while on the right limit point $b$ (if it exists) we get $\hat{\hat{f}}(b)=\lim_{y\to b-}f(y)$. So we conclude that the Legendre transform is an involution on the set of convex functions on $\mathbb{R}$ with those continuity properties.

We now come to the partial Legendre transform defined on a class of psh-functions with a certain symmetry, introduced by Kiselman in \cite{Kis}.

Let $U\subseteq \mathbb{C}^n$ be an open subset, and let $u(z,\tau)$ be a psh function on $U\times \mathbb{C}^*$ such that $u(z,\tau)=u(z,|\tau|)$. 

\begin{definition}
For $\lambda\in \mathbb{R}$ and $z\in U$ we define $$\hat{u}_{\lambda}(z):=\inf_{\tau\in \mathbb{C}^*}\{u(z,\tau)-\lambda\ln|\tau|^2\}.$$
\end{definition}

The family $\hat{u}_{\lambda}$ is known as the partial Legendre transform of $u$. To see why, let $y:=\ln|\tau|^2$. Since $u$ is psh and independent of the argument of $\tau$ we get that for a fixed $z\in U$: $f_z(y):=u(z,\tau)$ is convex in $y$, and we see that $\hat{u}_{\lambda}(z)=-\hat{f}_z(\lambda)$.

Thus $\hat{u}_{\lambda}$ is concave in $\lambda$, and by the involution property of the Legendre transform we get that

$$u(z,\tau)=\sup_{\lambda\in \mathbb{R}}\{\hat{u}_{\lambda}(z)+\lambda\ln|\tau|^2\}.$$

A fundamental property of the partial Legendre transform which is much less obvious is that $\hat{u}_{\lambda}$ is in fact psh. This follows directly from Kiselman's minimum principle \cite{Kis}, as $u-\lambda\ln|\tau|^2$ is psh and independent of the argument of $\tau$.

Let $v_{\lambda}$, $\lambda\in [a,b]$ be a concave family of psh functions on $U$. If $$u(z,\tau):=\sup_{\lambda\in [a,b]}\{v_{\lambda}(z)+\lambda\ln|\tau|^2\}$$ is psh then we see from the involution property of the Legendre transform that $\hat{u}_{\lambda}=v_{\lambda}$ for $\lambda\in [a,b]$ and $\hat{u}_{\lambda}=-\infty$ for $\lambda \notin [a,b]$. A supremum of psh functions is psh as long as it is u.s.c, and it follows from the elementary lemma below that this will be the case here as long as a certain boundedness condition is met.  

\begin{lemma} \label{usc}
If $v_{\lambda}$, $\lambda\in [a,b]$, is a concave family of u.s.c. functions such that $v_{\lambda}$ is locally bounded for $\lambda\in [a',b']$ where $a'<b'$, then $u:=\sup_{\lambda\in [a,b]}\{v_{\lambda}\}$ is also u.s.c.
\end{lemma}

\begin{proof}
Pick $x$. In a neighbourhood of $x$ we have by assumption local bounds $$v_{(a'+b')/2}(y)-v_{a'}(y)\leq C_1$$ and $$v_{(a'+b')/2}(y)-v_{b'}(y)\leq C_2.$$ By concavity this implies that for $\lambda\in[(a'+b')/2,b]$: $$\frac{\partial v_{\lambda}(y)}{\partial \lambda}\leq \frac{2C_1}{b'-a'}$$ while for $\lambda\in[a,(a'+b')/2]$: $$\frac{\partial v_{\lambda}(y)}{\partial \lambda}\geq -\frac{2C_2}{b'-a'}.$$

Given $N$ let $\lambda_j,$ $j=0,...,N$ be equidistributed points on $[a,b]$. Using the bounds we then get 
\begin{eqnarray*}
\limsup_{y\to x}(\sup_{\lambda\in [a,b]}\{v_{\lambda}(y)\})\leq \limsup_{y\to x}(\sup_{j=0,...,N}\{v_{\lambda_j}(y)\})+O(1/N)\leq \\ \leq \sup_{j=0,...,N}\{v_{\lambda_j}(x)\}+O(1/N)\leq \sup_{\lambda\in [a,b]}\{v_{\lambda}(x)\}+O(1/N),
\end{eqnarray*}
where the second inequality used that each $v_{\lambda}$ is u.s.c.
\end{proof}

Critical to our paper will be the fact that the Monge-Amp\`ere measure of $u$ can be understood in terms of the Monge-Amp\`ere measures of $\hat{u}_{\lambda}$. 

\begin{theorem} \label{Thmlocal}
If $\hat{u}_{\lambda}=-\infty$ for $\lambda\notin [a,b]$ while for $\lambda\in [a,b]$ each $\hat{u}_{\lambda}$ is locally bounded, then
$$(\pi_U)_*MA(u)=(n+1)\int_{\lambda=a}^{b}MA(\hat{u}_{\lambda})d\lambda,$$ where $\pi_U$ denotes the projection of $U\times \mathbb{C}^*$ to $U$. In particular $$\int_{U\times \mathbb{C}^*}MA(u)=(n+1)\int_a^b\left(\int_U MA(\hat{u}_{\lambda})\right)d\lambda.$$
\end{theorem} 

\begin{proof}
By simple scaling we can assume that $[a,b]=[0,1]$. The special case when $\hat{u}_{\lambda}=(1-\lambda)\phi+\lambda \psi-\lambda^2$ with $\phi$ psh and smooth and $\psi$ psh and locally bounded is a special case of \cite[Lem. 3.2]{WN3}. By approximation it is still true when also $\phi$ is just locally bounded. The proof also works for $\hat{u}_{\lambda}=(1-\lambda)\phi+\lambda \psi-\delta \lambda^2,$ $\delta>0$. 

For general $u$ we let for a given $N\in \mathbb{}N$, $j=0,...,2^N$ and $\lambda\in [j/2^N,(j+1)/2^N]$:

$$v^N_{\lambda}:=(1-(2^N\lambda-j))\hat{u}_{j/2^N}+(2^N\lambda-j)\hat{u}_{(j+1)/2^N}-\lambda^2/N$$ We then get that $$u^N:=\sup_{\lambda\in [0,1]}\{v^N_{\lambda}+\lambda\ln|\tau|^2\}$$ increases almost everywhere to $u$ and so $MA(u^N)$ converges weakly to $MA(u)$. 

We also let $$u^N_j:=\sup_{\lambda\in [j/2^N,(j+1)/2^N]}\{v^N_{\lambda}+\ln|\tau|^2\},$$ so we know that $$(\pi_U)_*MA(u^N_j)=(n+1)\int_{\lambda=j/2^N}^{(j+1)/2^N}MA(v^N_{\lambda})d\lambda.$$ 

Let $A(a,b):=\{\frac{\partial u^N}{\partial y}\in (a/2^N,b/2^N)\}$. These sets are all open in the plurifine topology since $$A(a,b)=\{u^N>\sup_{\lambda\in [0,1]\setminus (a/2^N,b/2^N)}\{v^N_{\lambda}+\lambda\ln|\tau|^2\}\}.$$

Note that $u^N=u^N_j$ on $A(j,j+1)$, and thus $$\mathbbm{1}_{A(j,j+1)}MA(u^N_j)=\mathbbm{1}_{A(j,j+1)}MA(u^N).$$ On $A(-1,j)$ we have instead that $$u^N_j=v^N_{j/2N}+(j/2^N)\ln|\tau|^2,$$ and thus $$\mathbbm{1}_{A(-1,j)}MA(u^N_j)=\mathbbm{1}_{A(-1,j)}MA(v^N_{j/2N}+(j/2^N)\ln|\tau|^2)=0.$$ Similarly one sees that $$\mathbbm{1}_{A(j+1,2^N+1)}MA(u^N_j)=\mathbbm{1}_{A(j+1,2^N+1)}MA(\phi_{j+1})=0.$$

On $A(j-\epsilon,j+\epsilon)$ we have that 
\begin{eqnarray*}
u^N=\sup_{\lambda\in [(j-\epsilon)/2^N,(j+\epsilon)/2^N]}\{v^N_{\lambda}+\lambda\ln|\tau|^2\}=:u^N_{j,\epsilon},
\end{eqnarray*}
 and so $$\mathbbm{1_{A(j-\epsilon,j+\epsilon)}}MA(u^N)=\mathbbm{1_{A(j-\epsilon,j+\epsilon)}}MA(u^N_{j,\epsilon}).$$ It can be easily showed that for any relatively compact set $K\subseteq U$ one can find a constant $C$ such that $\int_{K\times \mathbb{C}^*}MA(u^N_{j,\epsilon})\leq C\epsilon$. Thus $MA(u^N)$ puts no mass on the sets $\{\frac{\partial u^N}{\partial y}=j/2^N)\}$ and the same kind of argument shows that the same is true for $MA(u^N_j)$.

This lets us conclude that $$MA(u^N)=\sum_{j=0}^{2^N-1} MA(u^N_j)$$ and so $$(\pi_U)_*MA(u^N)=(n+1)\int_{\lambda=0}^{1}MA(v^N_{\lambda})d\lambda.$$ Finally, by letting $N\to \infty$ we get the Theorem.  
\end{proof}

\begin{remark}
In the forthcoming paper \cite{BWN} it is shown that more generally $$(\pi_U)_*(dd^cu)^{k+1}=(k+1)\int_{\lambda=a}^{b}(dd^c\hat{u}_{\lambda})^kd\lambda.$$ The proof in \cite{BWN}, rather than relying on \cite{WN3}, uses an explicit formula for $dd^c\hat{u}_{\lambda}$ essentially due to Kiselman \cite{Kis}, and is more direct.
\end{remark}

If $\theta$ is a smooth real $(1,1)$-form on a complex manifold $X$ and $u$ is $\pi_X^*\theta$-psh on $X\times \mathbb{C}^*$ then the associated Legendre transform becomes $$\hat{u}_{\lambda}(z):=\inf_{\tau\in \mathbb{C}^*}\{u(z,\tau)-\lambda\ln|\tau|^2\}.$$ It is easy to see, e.g. using a local potential for $\theta$, that everything we established for psh-functions works equally well for $\theta$-psh functions, e.g. that $\hat{u}_{\lambda}$ is a concave family of $\theta$-psh functions and $$u(z,\tau)=\sup_{\lambda\in \mathbb{R}}\{\hat{u}_{\lambda}(z)+\lambda\ln|\tau|^2\}.$$ The decomposition of the Monge-Amp\`ere of $u$ in terms of the Monge-Amp\`ere of $\hat{u}_{\lambda}$ will also look the same. 

\begin{theorem} \label{Thmtrans}
If $\hat{u}_{\lambda}=-\infty$ for $\lambda\notin [a,b]$ while for $\lambda\in [a,b]$ each $\hat{u}_{\lambda}$ is locally bounded away from a proper analytic subset $A\subseteq X$, then
$$(\pi_X)_*MA_{\pi_X^*\theta}(u)=(n+1)\int_{\lambda=a}^{b}MA_{\theta}(\hat{u}_{\lambda})d\lambda,$$ and in particular $$\int_{X\times \mathbb{C}^*}MA_{\pi_X^*\theta}(u)=(n+1)\int_a^b\left(\int_X MA_{\theta}(\hat{u}_{\lambda})\right)d\lambda.$$
\end{theorem}

This follows from Theorem \ref{Thmlocal} together with the fact that the Monge-Amp\`ere measures put no mass on the analytic subsets $A$ and $A\times \mathbb{C}^*$.

\begin{remark}
More generally one can consider the case when $u$ is psh (or $\theta$-psh) on $U\times (\mathbb{C}^*)^k$ and $u(z,\tau_1,...,\tau_k)=u(z,|\tau_1|,...,|\tau_k|)$. Given $\lambda \in \mathbb{R}^k$ one defines $$\hat{u}_{\lambda}(z):=\inf_{\tau\in \mathbb{C}^*}\{u(z,\tau_1,...,\tau_k)-\sum_i \lambda_i\ln|\tau_i|^k\}.$$ All the results above generalize to this setting.
\end{remark}

\section{Step 1: Volume formulas} \label{Sec:form}

Let $\alpha$ be a big $(1,1)$-class in $H^{1,1}(X,\mathbb{R})$, and let $Y$ be a smooth hypersurface that intersects $K(\alpha).$ Let also $\mu: \mathcal{X} \to X\times \mathbb{P}^1$ be the (standard) deformation of $X$ to the normal bundle of $Y$, as described in the introduction.

Pick positive constants $b,c$ such that $\alpha-c\{Y\}$ is big and $b>c$, and let $$\beta:=(\pi_X\circ\mu)^*\alpha+b\{X_0\}-c\{\mathcal{E}\}.$$ 

It is clear that $\beta\in H^{1,1}(\mathcal{X},\mathbb{R})$ restricts to $\alpha$ on $X_{\tau}$, $\tau\neq 0$.  

\begin{proposition} \label{Prop:kl}
The class $\beta$ is big, and $K(\beta)$ intersects each fiber $X_{\tau}$, $\tau\neq 0$, as well as $X'$ and $\mathcal{E}$.
\end{proposition}

\begin{proof}
To see that $\beta$ is big we note that $\beta$ also can be written $$\beta=\mu^*(\pi_X^*\alpha+(b-c)\{X_0\})+c\{X'\}.$$ The class $\pi_X^*\alpha+(b-c)\{X_0\}$ is easily seen to be big on $X\times \mathbb{P}^1$ and hence the pullback $\mu^*(\pi_X^*\alpha+(b-c)\{X_0\})$ is big on $\mathcal{X}$. The class $c\{X'\}$ is pseudoeffective, and since the sum of a big class with a pseudoeffective class is big, this shows that $\beta$ is big.

Let $\epsilon>0$ be small enough so that $\alpha-c\{Y\}-\epsilon[\omega]$ still is big, and let $\delta>0$ be small enough so that  class $$\gamma:=\epsilon((\pi_X\circ\mu)^*[\omega]+\{X_0\})-\delta\{\mathcal{E}\}$$ is K\"ahler. Then if $T$ is a closed positive current in $\beta-\gamma$  analytic singularities, the set where $T$ is smooth is contained in $K(\beta)$. Now note that 
\begin{eqnarray*}
\beta-\gamma=(\pi_X\circ\mu)^*(\alpha-\epsilon[\omega])+(b-\epsilon)\{X_0\}-(c-\delta)\{\mathcal{E}\}=\\=(\pi_X\circ\mu)^*(\alpha-\epsilon[\omega])+(b+\delta-c-\epsilon)\{X_0\}+(c-\delta)\{X'\}.
\end{eqnarray*}
Pick a closed positive current $S\in \alpha-\epsilon[\omega]$ with analytic singularities, and let $U\subseteq X$ denote the set where it is smooth. The closed positive current $$T:=(\pi_X\circ\mu)^*S+(b+\delta-c-\epsilon)(\pi_{\mathbb{P}^1}\circ \mu)^*\omega_{FS}+(c-\delta)[X']$$ then lies in $\beta-\gamma$ and has analytic singularities, and it is smooth on $\mu^{-1}(U\times \mathbb{P}^1)\setminus X'$. From this we see that $K(\beta)$ intersects each fiber $X_{\tau}$, $\tau\neq 0$, as well as $\mathcal{E}$. 

Let now $S'$ be a closed positive current in $\alpha-(c-\delta)\{Y\}-\epsilon[\omega]$ with analytic singularities, and let $U'$ be the set where it is smooth. The closed positive current 
\begin{eqnarray*}
T':=(\pi_X\circ\mu)^*(S'+(c-\delta)[Y])+(b-\epsilon)(\pi_{\mathbb{P}^1}\circ \mu)^*\omega_{FS}-(c-\delta)[\mathcal{E}]=\\=(\pi_X\circ\mu)^*S'+(b-\epsilon)(\pi_{\mathbb{P}^1}\circ \mu)^*\omega_{FS}
\end{eqnarray*} 
then lies in $\beta-\gamma$ and has analytic singularities, and it is smooth on $\mu^{-1}((U'\setminus Y)\times \mathbb{P}^1)$, which shows that $K(\beta)$ intersects $X'$.
\end{proof}

The goal of this section is to prove Theorem \ref{volcalc2}, i.e. the four formulas:

\begin{enumerate}
\item $$\textrm{vol}(\beta)=(n+1)\left((b-c)\textrm{vol}(\alpha)+\int_0^c\textrm{vol}(\alpha-t\{Y\})dt\right),$$
\item $$\textrm{vol}_{\mathcal{X}|X_{\tau}}(\beta)=\textrm{vol}(\alpha), \qquad{} \forall \tau\neq 0,$$
\item $$\textrm{vol}_{\mathcal{X}|X'}(\beta)=\textrm{vol}(\alpha-c\{Y\}),$$
\item $$\textrm{vol}_{\mathcal{X}|\mathcal{E}}(\beta)=n\int_0^c\textrm{vol}_{X|Y}(\alpha-t\{Y\})dt.$$
\end{enumerate}

Let $\theta$ be a smooth form representing $\alpha$. Let $s_{\mathcal{E}}$ be a defining section for $\mathcal{E}$ and let $h$ be a smooth $S^1$-invariant hermitian metric on the associated line bundle, which we for simplicity choose so that $|s_{\mathcal{E}}|_h\leq 1$. Note that $$dd^c\ln|s_{\mathcal{E}}|^2_h=[\mathcal{E}]-\eta,$$ where $\eta$ is the curvature form of $h$. Let $\omega_{FS}=dd^c\ln(1+|\tau|^2)$ denote the Fubini-Study form on the base $\mathbb{P}^1$. Thus $$\tilde{\theta}:=(\pi_X\circ \mu)^*\theta+b(\pi_{\mathbb{P}^1}\circ \mu)^*\omega_{FS}-c\eta$$ is a smooth $S^1$-invariant form representing $\beta$.

Let $$\tilde{u}:=\sup\{v\leq 0: \tilde{\psi}\in PSH(\mathcal{X},\tilde{\theta})\}.$$ Clearly $\tilde{u}$ has minimal singularities and hence its Monge-Amp\`ere can be used to calculate the volume and restricted volumes of $\beta$.

When performing the partial Legendre transform it is easier to work on $(X\times \mathbb{C}^*,\pi_X^*\theta)$ than on $(\mathcal{X},\tilde{\theta})$. Thus we introduce a similar looking envelope but now in $PSH(X\times \mathbb{C}^*,\pi_X^*\theta)$: $$u:=\sup\{v\leq b\ln(1+|\tau|^2)+c\ln|s_{\mathcal{E}}|^2_h: \psi\in PSH(X\times \mathbb{C}^*, \pi_X^*\theta)\}.$$ 

Note that since $b\ln(1+|\tau|^2)+c\ln|s_{\mathcal{E}}|^2_h$ and $\pi_X^*\theta$ is $S^1$-invariant (i.e. only depends on $|\tau|$) $u$ is also $S^1$-invariant.

$u$ and $\tilde{u}$ are related in a simple way:

\begin{proposition} \label{proptilde}
$$\tilde{u}=u-b\ln(1+|\tau|^2)-c\ln|s_{\mathcal{E}}|^2_h$$ on $X\times \mathbb{C}^*$.
\end{proposition}

\begin{proof}
Clearly $\tilde{u}+b\ln(1+|\tau|^2)+c\ln|s_{\mathcal{E}}|^2_h\leq b\ln(1+|\tau|^2)+c\ln|s_{\mathcal{E}}|^2_h$ and it lies in $PSH(X\times \mathbb{C}^*,\pi_X^*\theta)$, so by definition we get that $$\tilde{u}+b\ln(1+|\tau|^2)+c\ln|s_{\mathcal{E}}|^2_h \leq u.$$ On the other hand $$u-b\ln(1+|\tau|^2)-c\ln|s_{\mathcal{E}}|^2_h\leq 0$$ and it lies in $PSH(X\times \mathbb{C}^*,\tilde{\theta}\}$. Being bounded it extends to a $\tilde{\theta}$-psh function on the whole of $\mathcal{X}$, and so by definition $$u-b\ln(1+|\tau|^2)-c\ln|s_{\mathcal{E}}|^2_h\leq \tilde{u}.$$ 

\end{proof}

Let us consider the partial Legendre transform $\hat{u}_{\lambda}$ of $u$, i.e. $$\hat{u}_{\lambda}(z):=\inf_{|\tau|>0}\{u(z,\tau)-\lambda\ln|\tau|^2\}.$$ Recall from Section \ref{Sec:trans} that $\hat{u}_{\lambda}\in PSH(X,\theta)$. Since $c\ln|s_{\mathcal{E}}|^2_h\leq 0$ we get that $u\leq b\ln(1+|\tau|^2)$ and so $$\hat{u}_{\lambda}\leq \inf_{|\tau|>0}\{b\ln(1+|\tau|^2)+\lambda \ln|\tau|^2\}=-\infty$$ for $\lambda<0$ or $\lambda>b$ thus $u_{\lambda}=-\infty$ when  $\lambda\notin[0,b]$. By the involution property we thus have that $$u(z,\tau)=\sup_{\lambda\in [0,b]}\{\hat{u}_{\lambda}(z)+\lambda\ln|\tau|^2\}.$$

We now want to analyse $\hat{u}_{\lambda}$. To do this we will compare $\hat{u}_{\lambda}$ with related envelopes $v_t$, to be introduced shortly. 

First note that $\ln|s_{\mathcal{E}}|^2_h$ descends to $X\times \mathbb{C}$ since it is $-\infty$ on $\mathcal{E}$. Restricted to $X\times\{0\}\equiv X'$ $s_{\mathcal{E}}$ is a defining section for $Y$, so we write the restriction of $\ln|s_{\mathcal{E}}|^2_h$ to $X\times\{0\}\equiv X'$ as $\ln|s_Y|^2_h$. We then have that $$dd^c\ln|s_Y|^2_h=[Y]-\sigma$$ where $\sigma$ denotes the restriction of $\eta$ to $X\equiv X'$.

\begin{definition}
For $t\in [0,c]$: $$v_t:=\sup\{v\leq t\ln|s_Y|^2_h: v\in PSH(X,\theta)\}.$$
\end{definition}

Note that since $\ln|s_Y|^2_h\leq 0$, $v_t$ is decreasing in $t$.

Clearly $v_t$ has Lelong number at least $t$ along $Y$, and among such $\theta$-psh functions it has minimal singularities. 

\begin{proposition} We have that $v_t$ is locally bounded on $K(\alpha-t\{Y\})\setminus Y$,
$$\int_X MA_{\theta}(v_t)=\textrm{vol}(\alpha-t\{Y\})$$ and $$\int_Y MA_{(\theta-t\sigma)_{|Y}}((v_t-t\ln|s_Y|^2_h)_{|Y})=\textrm{vol}_{X|Y}(\alpha-t\{Y\}).$$
\end{proposition}

\begin{proof}
$\theta-t\sigma$ is a smooth form representing $\alpha-\{Y\}$. Let $$\tilde{v}:=\sup\{v\leq 0: \in PSH(X,\theta-t\sigma)\}.$$ Then $\tilde{v}$ has minimal singularities so $$\int_X MA_{\theta-t\sigma}(\tilde{v})=\textrm{vol}(\alpha-t\{Y\}).$$ Since $\tilde{v}+t\ln|s_Y|^2_h\leq t\ln|s_Y|^2_h$ is $\theta$-psh we get that $\tilde{v}+t\ln|s_Y|^2_h\leq \psi_t$. On the other hand, $v_t-\lambda\ln|s_Y|^2_h\leq 0$ and is $\theta-t\sigma$-psh, which shows that $v_t=\tilde{v}+\lambda\ln|s_Y|^2_h$. It follows that $v_t$ is locally bounded on $K(\alpha-t\{Y\})\setminus Y$, $$MA_{\theta}(v_t)=MA_{\theta-t\sigma}(\tilde{v})$$ and $$MA_{(\theta-t\sigma)_{|Y}}((v_t-t\ln|s_Y|^2_h)_{|Y})=MA_{(\theta-t\sigma)_{|Y}}(\tilde{v}_{|Y})$$ and so the Proposition follows. 
\end{proof}

\begin{proposition} \label{equivalence}
For $\lambda\in [0,b]$ and $t:=\max(c-\lambda,0)$ one can find a constant $C$ so that $$|\hat{u}_{\lambda}-v_t|\leq C.$$ 
\end{proposition}

\begin{proof}
Note that 
\begin{equation} \label{eq1}
\hat{u}_{\lambda}(z)+\lambda\ln|\tau|^2\leq u(z,\tau)\leq b\ln(1+|\tau|^2)+c\ln|s_{\mathcal{E}}|^2_h.
\end{equation} 
on $X\times \mathbb{C}^*$ and since the right-hand-side is bounded near $X\times \{0\}$ the inequality extends to $X\times \mathbb{C}$. 
The Lelong number along $Y\times \{0\}$ of the right hand side is $c$ and by the monotonicity of Lelong numbers the Lelong number along $Y\times \{0\}$ of the left hand side must be bigger than or equal to $c$. On the other hand is it easy to see that this Lelong number equals $\nu_Y(\hat{u}_{\lambda})+\lambda,$ giving us that $$\nu_Y(\hat{u}_{\lambda})\geq t.$$ From this it follows that $\hat{u}_{\lambda}\leq t\ln|s_Y|^2_h+C$ for some constant $C$, and thus $\hat{u}_{\lambda}\leq v_t+C$. 

Also note that $v_t+\lambda\ln|\tau|^2\in PSH(X\times \mathbb{C})$ and that $$\nu_{Y\times \{0\}}(v_t+\lambda\ln|\tau|^2)=\nu_Y(v_t)+\lambda\geq t+\lambda\geq c.$$ It follows that $v_t+\lambda\ln|\tau|^2\leq c\ln|s_{\mathcal{E}}|^2_h+C'$ when $|\tau|\leq 1$ say, for some constant $C'$. Since $c\ln|s_{\mathcal{E}}|^2_h$ is bounded from below by some constant for $|\tau|>1$ we get that $$v_t+\lambda\ln|\tau|^2\leq b\ln(1+|\tau|^2)+c\ln|s_{\mathcal{E}}|^2_h+C''$$ on $X\times \mathbb{C}$ for some constant $C''$. It follows that $$v_t+\lambda\ln|\tau|^2\leq u+C''$$ which implies that $$u_t\leq \inf_{|\tau|>0}\{u-\lambda\ln|\tau|^2\}+C''=\hat{u}_{\lambda}+C'',$$ and so the Proposition follows.
\end{proof}

\begin{corollary} \label{corgreat}
For $\lambda\in [0,b]$ and $t:=\max(c-\lambda,0)$ we have that $\hat{u}_{\lambda}$ is locally bounded on $K(\alpha-t\{Y\})\setminus Y$,
$$\int_X MA_{\theta}(\hat{u}_{\lambda})=\textrm{vol}(\alpha-t\{Y\})$$ and $$\int_Y MA_{(\theta-t\sigma)_{|Y}}((\hat{u}_{\lambda}-t\ln|s_Y|^2_h)_{|Y})=\textrm{vol}_{X|Y}(\alpha-t\{Y\}).$$ 
\end{corollary}

Since $v_t$ was seen to be decreasing in $t$, another consequence of Proposition \ref{equivalence} is that $\hat{u}_{\lambda}$ is increasing in $\lambda$ up to constants, i.e. that if $\lambda_1\leq \lambda_2$ then for some constant $C$ we have that $$\hat{u}_{\lambda_1}\leq \hat{u}_{\lambda_2}+C.$$ In particular each $\hat{u}_{\lambda}$ is locally bounded away from the proper analytic subset $E_{nK}(\alpha-c\{Y\})\cup Y$.

We are now ready to prove the first formula of Theorem \ref{volcalc2}, i.e. that $$\textrm{vol}(\beta)=(n+1)\left((b-c)\textrm{vol}(\alpha)+\int_0^c\textrm{vol}(\alpha-t\{Y\})dt\right).$$

\begin{proof}[Proof of Theorem \ref{volcalc2}, formula 1]
Using Proposition \ref{proptilde} we see that $$MA_{\tilde{\theta}}(\tilde{u})=MA_{\pi_X^*\theta}(u)$$ on $X\times \mathbb{C}^*$. Since $MA_{\tilde{\theta}}(\tilde{u})$ puts no mass on $X_0$ and $X_{\infty}$ we get that $$\textrm{vol}(\beta)=\int_{\mathcal{X}}MA_{\tilde{\theta}}(\tilde{u})=\int_{X\times \mathbb{C}^*}MA_{\pi_X^*\theta}(u).$$ On the other hand, combining Theorem \ref{Thmtrans} and Corollary \ref{corgreat} yields 
 
\begin{eqnarray*}
\int_{X\times \mathbb{C}^*}MA_{\pi_X^*\theta}(u)=(n+1)\int_{\lambda=0}^b\left(\int_X MA_{\theta}(\hat{u}_{\lambda})\right)d\lambda=\\=(n+1)\left((b-c)\textrm{vol}(\alpha)+\int_0^c\textrm{vol}(\alpha-t\{Y\})dt\right).
\end{eqnarray*}
\end{proof}

If we apply the formula on $\beta$ and $\beta-\{X_0\}$ we get the following corollary, which we will have use for later.

\begin{corollary} \label{useful}
$$\textrm{vol}(\beta)-\textrm{vol}(\beta-\{X_0\})=(n+1)\textrm{vol}(\alpha).$$
\end{corollary}

We will move on to the second formula of Theorem \ref{volcalc2}, i.e. the claim that for $\tau\neq 0$ we have $$\textrm{vol}_{\mathcal{X}|X_{\tau}}(\beta)=\textrm{vol}(\alpha).$$

\begin{proof}[Proof of Theorem \ref{volcalc2}, formula 2]
Since $\tilde{u}\in PSH(\mathcal{X},\tilde{\theta})$ has minimal singularities we know that $$\textrm{vol}_{\mathcal{X}|X_{\tau}}(\beta)=\int_{X_{\tau}}MA_{\tilde{\theta}_{|X_{\tau}}}(\tilde{u}_{|X_{\tau}})=\int_X MA_{\theta}(u_{|X_{\tau}}).$$ Using Proposition \ref{equivalence} and the fact that $u(z,\tau)\geq \hat{u}_b(z)+b\ln|\tau|^2$ we see that $$u_{|X_{\tau}}\geq \hat{u}_b+b\ln|\tau|^2\geq v_0-C$$ for some constant $C$. Since $v_0\in PSH(X,\theta)$ has minimal singularities it follows that $u_{|X_{\tau}}$ also has minimal singularities, and hence $$\int_X MA_{\theta}(u_{|X_{\tau}})=\textrm{vol}(\alpha).$$
\end{proof}

We have seen that on $X\times \mathbb{C}^*$, $\phi$ and hence also $\tilde{u}$ could be expressed in terms of the partial Legendre transform $\hat{u}_{\lambda}$. We now want a similar description of $\tilde{u}$ restricted to $X'$ and $\mathcal{E}$ in order to calculate the restricted volume of $\beta$ along those hypersurfaces. 

We start by looking at what happens on $X'$.

Pick a point $z\in K(\alpha-c\{Y\})\setminus Y$ and let us identify $(z,0)$ with the corresponding point in $X'$. Since $\tilde{u}$ is quasi-psh we have that $$\tilde{u}(z,0)=\lim_{\epsilon\to 0}\max_{|\tau|=\epsilon}\{\tilde{u}(z,\tau)\}=\lim_{\tau\to 0}\tilde{u}(z,\tau),$$ where the last equality comes from the fact that $\tilde{u}$ is $S^1$-invariant.

Note that $z\in K(\alpha-c\{Y\})$ means that $\hat{u}_0(z)>-\infty$. Also note that we have the uniform bound $\hat{u}_{\lambda}(z)\leq b\ln 2.$ It follows that for any $\delta>0$ and $|\tau|<1$ small enough
\begin{eqnarray*}
\hat{u}_0(z)\leq \sup_{\lambda\in [0,b]}\{\hat{u}_{\lambda}(z)+\lambda\ln|\tau|^2\}\leq \sup_{\lambda\in [0,\delta]}\{\hat{u}_{\lambda}(z)\}.
\end{eqnarray*}
From this we get that 
\begin{eqnarray*}
\hat{u}_0(z)-c\ln|s_Y|^2_h(z) \leq \tilde{u}(z,0) \leq \sup_{\lambda\in [0,\delta]}\{\hat{u}_{\lambda}(z)\}-c\ln|s_Y|^2_h(z).
\end{eqnarray*}
To show that in fact 
\begin{eqnarray} \label{id1}
\tilde{u}(z,0)=\hat{u}_0(z)-c\ln|s_Y|^2_h(z)
\end{eqnarray}
we just need the following lemma.

\begin{lemma}
$$\hat{u}_0=\lim_{\delta\to 0}\sup_{\lambda\in [0,\delta]}\{\hat{u}_{\lambda}\}.$$
\end{lemma}

\begin{proof}
We note that $\phi_{\delta}:=\sup_{\lambda\in [0,\delta]}^*\{\hat{u}_{\lambda}+\lambda\ln|\tau|^2\}$ is $\pi_X^*\theta$-psh and bounded from above by $\phi$. Note that 
\begin{eqnarray} \label{eqest1}
|\phi_{\delta}(z,\tau)-\sup_{\lambda\in [0,\delta]}^*\{\hat{u}_{\lambda}(z)\}|\leq \delta|\ln|\delta|^2|.
\end{eqnarray}
As $\delta\to 0$, $\phi_{\delta}$ decreases to a $\pi_X^*\theta$-psh function $\phi_0\leq u$, and by (\ref{eqest1}) $$\phi_0(z,\tau)=\lim_{\delta\to 0}\sup_{\lambda\in [0,\delta]}^*\{\hat{u}_{\lambda}(z)\}.$$ Thus $$\lim_{\delta\to 0}\sup_{\lambda\in [0,\delta]}\{\hat{u}_{\lambda}(z)\}\leq \phi_0(\tau,z)$$ and taking $\inf$ over $|\tau|>0$ gives us that 
$$\lim_{\delta\to 0}\sup_{\lambda\in [0,\delta]}\{\hat{u}_{\lambda}(z)\}\leq \inf_{|\tau|>0}u(z,\tau)\leq \inf_{|\tau|>0}u(z,\tau)=\hat{u}_0(z).$$
\end{proof}

Given that we showed (\ref{id1}) on $K(\alpha-c\{Y\})\setminus Y$ whose complement in $X'=X$ is contained in a proper analytic set the identity (\ref{id1}) extends to the whole of $X'$.

We are now ready to prove the third formula of Theorem \ref{volcalc2}, i.e. that
$$\textrm{vol}_{\mathcal{X}|X'}(\beta)=\textrm{vol}(\alpha-c\{Y\}).$$

\begin{proof}[Proof of Theorem \ref{volcalc2}, formula 3]
We have that $$\textrm{vol}_{\mathcal{X}|X'}(\beta)=\int_{X'}MA_{\tilde{\theta}_{|X'}}(\tilde{u}_{|X'}).$$ Given the identification $X'=X$ we have that $\tilde{\theta}_{|X'}=\theta-c\eta'$ and $\tilde{u}_{|X'}=\hat{u}_0-c\ln|s_Y|^2_h$, which means that $$MA_{\tilde{\theta}_{|X'}}(\tilde{u}_{|X'})=MA_{\theta}(\hat{u}_0)$$ and hence $$\textrm{vol}_{\mathcal{X}|X'}(\beta)=\int_X MA_{\theta}(\hat{u}_0)=\textrm{vol}(\alpha-c\{Y\})$$ by Corollary \ref{corgreat}.
\end{proof}

We now move on to analyse $\tilde{u}$ on $\mathcal{E}$.

Let $y$ be a point in $Y$ and pick local holomorphic coordinates $z_1,...,z_n$ centered at $Y$ so that $z_1=0$ is a local equation for $Y$. In this neighbourhood of $y$ in $X$ we have that $\ln|s_Y|^2_h=\ln|z_1|^2-g_Y$ where $dd^cg_Y=\eta'.$ 

We also have that $\tilde{z}_1:=z_1/\tau, z_2, ... , z_n, \tau$ are local holomorphic coordinates on $\mathcal{X}$ in a neighbourhood of the fibers of $\mathcal{E}$ over the points of $Y$ with coordinates $(0,z_2,...,z_n)$. In these coordinates $\mathcal{E}$ has defining equation $\tau=0$, and $(\tilde{z}_1,z_2,...,z_n,0)$ is the coordinate for the point in $\mathcal{E}\subseteq \mathcal{X}$ which is the limit of the points $(\tau \tilde{z}_1,z_2,...,z_n,\tau)\in X\times \mathbb{C}^*\subseteq \mathcal{X}$. Also note that in these coordinates $\ln|s_{\mathcal{E}}|^2_h=\ln|\tau|^2-g_{\mathcal{E}}$ where $dd^cg_{\mathcal{E}}=\eta$.

Recall that $$\tilde{u}=\sup_{\lambda\in [0,b]}\{\hat{u}_\lambda+\lambda\ln|\tau|^2\}-b\ln(1+|\tau|^2)-c\ln|s_{\mathcal{E}}|^2_h$$ on $X\times \mathbb{C}^*$.

In local coordinates we see that 
\begin{eqnarray*}
\phi_{\lambda}:=\hat{u}_{\lambda}+\lambda\ln|\tau|^2-b\ln(1+|\tau|^2)-c\ln|s_{\mathcal{E}}|^2_h=\\=\hat{u}_{\lambda}-(c-\lambda)\ln|s_Y|_h^2+(c-\lambda)\ln|\tilde{z}_1|^2-b\ln(1+|\tau|^2)-(c-\lambda)g_Y+cg_{\mathcal{E}}.
\end{eqnarray*}
First note that 
$$\phi_{\lambda}\leq u(z,1)+(\lambda-c)\ln|\tau|^2-b\ln(1+|\tau|^2)+cg_{\mathcal{E}}.$$ This shows that for $\lambda>c$: $\phi_{\lambda}$ tends to $-\infty$ as $\tau\to 0$, hence it extends as a quasi-psh function over $\tau=0$ with $\phi_{\lambda}(\tilde{z_1},z_2,...,z_n,0)=-\infty$. 

Now we look at the case $\lambda\in [0,c]$. Since $\hat{u}_{\lambda}$ has Lelong number at least $c-\lambda$ along $Y$ this function also extends as a quasi-psh function for $\tau=0$. Let $$\psi_{\lambda}:=(\hat{u}_{\lambda}-(c-\lambda)(\ln|s_Y|^2_h+g_Y))_{|Y}$$ then $$\phi_{\lambda}(\tilde{z_1},z_2,...,z_n,0)=\psi_{\lambda}\circ \pi_Y+(c-\lambda)\ln|\tilde{z}_1|^2+cg_{\mathcal{E}},$$ where $\pi_Y$ denotes the projection from $\mathcal{E}$ to $Y$.

Note that if $(z_2,...,z_n)\in K(\alpha-c\{Y\})\cap Y$ and $|\tilde{z}_1|>0$ then for $\lambda\in [0,c]$ $v_{\lambda}$ is locally bounded in a neighbourhood of $(\tilde{z}_1,z_2,...,z_n,0)$. Thus by Lemma \ref{usc} we get that $\sup_{\lambda\in [0,b]}\phi_{\lambda}$ is u.s.c. and then also quasi-psh in a neighbourhood of $(\tilde{z}_1,z_2,...,z_n,0).$ This means that the equality $$\tilde{u}=\sup_{\lambda\in [0,b]}\phi_{\lambda}$$ extends from $X\times \mathbb{C}^*$ to $\mathcal{E}\setminus A$, where $A$ is the proper analytic subset consisting of the zero and infinity section of $\mathcal{E}$ and $\pi_Y^{-1}(E_{nK}(\alpha-c\{Y\})\cap Y)$. 

To conclude, we have that on $\mathcal{E}\setminus A$, in local coordinates
\begin{eqnarray} \label{identity}
\tilde{u}=\sup_{\lambda\in[0,c]}\{\psi_{\lambda}\circ \pi_Y+(c-\lambda)\ln|\tilde{z}_1|^2\}+cg_{\mathcal{E}}.
\end{eqnarray} 

We are now ready to prove the last formula of Theorem \ref{volcalc2}, i.e. that $$\textrm{vol}_{\mathcal{X}|\mathcal{E}}(\beta)=n\int_0^c \textrm{vol}_{X|Y}(\alpha-t\{Y\})dt.$$

\begin{proof}[Proof of Theorem \ref{volcalc2}, formula 4]

We note that $$\tilde{\theta}_{|\mathcal{E}}=(\pi_Y)^*\theta_{|Y}-cdd^cg_{\mathcal{E}}$$ and thus $$(dd^c\tilde{u}+\tilde{\theta})_{|\mathcal{E}}=dd^c\sup_{\lambda\in[0,c]}\{\psi_{\lambda}\circ \pi_Y+(c-\lambda)\ln|\tilde{z}_1|^2\}+(\pi_Y)^*\theta_{|Y}.$$ We also note that $$dd^c \psi_{\lambda}+\theta_{|Y}=dd^c(\hat{u}_{\lambda}-(c-\lambda)\ln|s_Y|^2)+(\theta-(c-\lambda)\sigma)_{|Y}$$ and so by Corollary \ref{corgreat} $$\int_Y MA_{\theta_{|Y}}(\psi_{\lambda})=\textrm{vol}_{X|Y}(\alpha-(c-\lambda)\{Y\}).$$

Using Theorem \ref{Thmtrans} then finally gives us that 

\begin{eqnarray*}
\textrm{vol}_{\mathcal{X}|\mathcal{E}}(\beta)=\int_{\mathcal{E}}MA_{\tilde{\theta}_{|\mathcal{E}}}(\tilde{u}_{|\mathcal{E}})=\\=n\int_{\lambda=0}^c\left(\int_Y MA_{\theta_{|Y}}(u_{\lambda})\right)d\lambda=n\int_0^c \textrm{vol}_{X|Y}(\alpha-t\{Y\})dt.
\end{eqnarray*}

\end{proof}

\section{Step 2: The key equality}

The equality left to prove is $$\textrm{vol}_{\mathcal{X'}|X_0}(\beta)+\textrm{vol}_{\mathcal{X}|\mathcal{E}}(\beta)=\textrm{vol}(\alpha),$$ as this, given Theorem \ref{volcalc2}, will imply Theorem C. 

We start by showing the easy inequality $$\textrm{vol}_{\mathcal{X'}|X_0}(\beta)+\textrm{vol}_{\mathcal{X}|\mathcal{E}}(\beta)\leq \textrm{vol}(\alpha).$$

\begin{proof}[Proof of Theorem \ref{keyeq}, easy inequality]

Pick $u\in PSH(\mathcal{X},\tilde{\theta})$ with analytic singularities. Locally, away from the singularities of $u$, the measures $MA_{\tilde{\theta}_{|X_\tau}}(u_{|X_{\tau}})$ converge weakly to the measure $MA_{\tilde{\theta}_{|X'}}(u_{|X'})$ on $X'$ and $MA_{\tilde{\theta}_{|\mathcal{E}}}(u_{|\mathcal{E}})$ on $\mathcal{E}$. This implies that $$ \liminf_{\tau\to 0} \int_{X_{\tau}}MA_{\tilde{\theta}_{|X_\tau}}(u_{|X_{\tau}})\geq \int_{X'}MA_{\tilde{\theta}_{|X'}}(u_{|X'})+\int_{\mathcal{E}}MA_{\tilde{\theta}_{|\mathcal{E}}}(u_{|\mathcal{E}}).$$ Since $$ \int_{X_{\tau}}MA_{\tilde{\theta}_{|X_\tau}}(u_{|X_{\tau}})\leq \textrm{vol}_{\mathcal{X}|X_{\tau}}(\beta)=\textrm{vol}(\alpha)$$ it follows that $$\textrm{vol}_{\mathcal{X'}|X_0}(\beta)+\textrm{vol}_{\mathcal{X}|\mathcal{E}}(\beta)\leq\textrm{vol}(\alpha).$$
\end{proof}

\section{Step 3: Measure control}

Recall from Section \ref{Sec:form} that $\tilde{\theta}$ is a smooth $S^1$-invariant form representing $\beta$.

Let $\chi: \mathbb{R} \to \mathbb{R}$ be a smooth convex function such that $\chi(x)=x$ for $x\geq 0$ and $\chi(x)=e^x-1$ for $x\leq -1$. 

We then let $g_k:=\chi(\ln|\tau|^2+k)-k-\ln(1+|\tau|^2)$. Note that $g_k$ is $\omega_{FS}$-psh, decreases to $\ln|\tau|^2-\ln(1+|\tau|^2)$, and that $g_k=\ln|\tau|^2-\ln(1+|\tau|^2)$ for $|\tau|>e^{-k/2}$.

Let $\omega:=dd^c\chi(\ln|\tau|^2)$ and $\omega_k:=\omega_{FS}+dd^cg_k$. We note that $\omega_k$ are probability measures on $\mathbb{P}^1$ and that 
$$\omega_k=(e^{-k/2})_*\omega$$ where $e^{-k/2}$ here denotes the function $\tau \mapsto e^{-k/2}\tau$.  

Now let $u_k$ be the associated envelopes $$u_k:=\sup\{u\leq g_k: \phi\in PSH(\mathcal{X},\tilde{\theta})\}.$$ 

Since $g_k$ decreases to $\ln|\tau|^2-\ln(1+|\tau|^2)$ it follows that $u_k$ decreases to some $u_{\infty}\in PSH(\mathcal{X}, \tilde{\theta})$ such that $$u_{\infty}\leq \ln|\tau|^2-\ln(1+|\tau|^2).$$ 
 
\begin{proposition} \label{Propineq}
$$\int_{\mathcal{X}}MA_{\tilde{\theta}}(u_{\infty})\leq \textrm{vol}(\beta-\{X_0\}).$$
\end{proposition}

\begin{proof}
As $u_{\infty}\leq \ln|\tau|^2-\ln(1+|\tau|^2)$ it follows that $v:=u_{\infty}-\ln|\tau|^2+\ln(1+|\tau|^2)$ is $\tilde{\theta}-(\pi_{\mathbb{P}^1}\circ \mu)^*\omega_{FS}$-psh. Note that $\tilde{\theta}-(\pi_{\mathbb{P}^1}\circ \mu)^*\omega_{FS} \in \beta-\{X_0\}$ and that $MA_{\tilde{\theta}}(u_{\infty})=MA_{\tilde{\theta}-(\pi_{\mathbb{P}^1}\circ \mu)^*\omega_{FS}}(v)$. It follows that $$\int_{\mathcal{X}}MA_{\tilde{\theta}}(u_{\infty})=\int_{\mathcal{X}}MA_{\tilde{\theta}-(\pi_{\mathbb{P}^1}\circ \mu)^*\omega_{FS}}(v)\leq \textrm{vol}(\beta-\{X_0\}).$$
\end{proof}

Let $D_k:=\{u_k=g_k\}$ denote the contact set and recall that by Theorem \ref{DNT} $$MA_{\tilde{\theta}}(u_k)=\mathbbm{1}_{D_k}(\tilde{\theta}+dd^c g_k)^{n+1}.$$

\begin{proposition} For any $\epsilon>0$ we have that
$$\limsup_{k\to \infty}\int_{\{|\tau|>\epsilon\}}MA_{\tilde{\theta}}(u_k)\leq \textrm{vol}(\beta-\{X_0\}).$$
\end{proposition}

\begin{proof}
Note that for $k$ large $g_k=\ln|\tau|^2-\ln(1+|\tau|^2)$ on $\{|\tau|>\epsilon\}$. Since $u_k$ is decreasing it follows that $D_k\cap \{|\tau|>\epsilon\}$ is decreasing to $D_{\infty}\cap\{|\tau|>\epsilon\}$ where $D_{\infty}:=\cap_{k=0}^{\infty}D_k$. We thus get that 
\begin{eqnarray*}
\mathbbm{1}_{\{|\tau|>\epsilon\}}MA_{\tilde{\theta}}(u_k)=\mathbbm{1}_{D_k\cap\{|\tau|>\epsilon\}}(\tilde{\theta}-(\pi_{\mathbb{P}^1}\circ \mu)^*\omega_{FS})^{n+1}\to \\ \to \mathbbm{1}_{D_{\infty}\cap\{|\tau|>\epsilon\}}(\tilde{\theta}-(\pi_{\mathbb{P}^1}\circ \mu)^*\omega_{FS})^{n+1}.
\end{eqnarray*}
On the other hand we know that $MA_{\tilde{\theta}}(u_k)$ converges weakly to $MA_{\tilde{\theta}}(u_{\infty})$ where $u_{\infty}$ is locally bounded, and it follows that for any constant $C$: $$\mathbbm{1}_{\{|\tau|>\epsilon, u_{\infty}>-C\}}MA_{\tilde{\theta}}(u_{\infty})=\mathbbm{1}_{D_{\infty}\cap\{|\tau|>\epsilon, \phi_{\infty}>-C\}}(\tilde{\theta}-(\pi_{\mathbb{P}^1}\circ \mu)^*\omega_{FS})^{n+1}$$ and hence letting $C\to \infty$: $$\mathbbm{1}_{\{|\tau|>\epsilon\}}MA_{\tilde{\theta}}(u_{\infty})=\mathbbm{1}_{D_{\infty}\cap\{|\tau|>\epsilon\}}(\tilde{\theta}-(\pi_{\mathbb{P}^1}\circ \mu)^*\omega_{FS})^{n+1}.$$ From this it follows that
\begin{eqnarray*}
\lim_{k\to \infty}\int_{\{|\tau|>\epsilon\}}MA_{\tilde{\theta}}(u_k)=\int_{\{|\tau|>\epsilon\}}MA_{\tilde{\theta}}(u_{\infty})\leq \textrm{vol}(\beta-\{X_0\}),
\end{eqnarray*}
where the last inequality comes from Proposition \ref{Propineq}

\end{proof}

\begin{proposition} \label{propliminf}
$$\liminf_{k\to \infty}\int_{\{|\tau|\leq \epsilon\}}MA_{\tilde{\theta}}(u_k)\geq (n+1)\textrm{vol}(\alpha).$$
\end{proposition}

\begin{proof}
We have that
\begin{eqnarray*}
\liminf_{k\to \infty}\int_{\{|\tau|\leq \epsilon\}}MA_{\tilde{\theta}}(u_k)=\liminf_{k\to \infty}\left(\int_{\mathcal{X}}MA_{\tilde{\theta}}(u_k)-\int_{\{|\tau|\leq \epsilon\}}MA_{\tilde{\theta}}(u_k)\right)=\\=\textrm{vol}(\beta)-\limsup_{k\to \infty}\int_{\{|\tau|>\epsilon\}}MA_{\tilde{\theta}}(u_k)\geq \textrm{vol}(\beta)-\textrm{vol}(\beta-\{X_0\})=(n+1)\textrm{vol}(\alpha),
\end{eqnarray*}
where in last equality we used Corollary \ref{useful}.
\end{proof}

\begin{proposition} \label{proplimit}
$$\liminf_{k\to \infty}\int_{\mathbb{P}^1}\left(\int_{X_{\tau}\cap D_k}\tilde{\theta}^n\right)\omega_k\geq \textrm{vol}(\alpha).$$
\end{proposition}

\begin{proof}
Note that $dd^cg_k\leq \omega_{FS}+dd^cg_k=\omega_k$ and that $$(\tilde{\theta}+\omega_k)^{n+1}=\tilde{\theta}^{n+1}+(n+1)\tilde{\theta}^n\wedge \omega_k.$$ Using this we get that
$$MA_{\tilde{\theta}}(u_k)=\mathbbm{1}_{D_k}(\tilde{\theta}+dd^cg_k)^{n+1}\leq \mathbbm{1}_{D_k}(\tilde{\theta}+\omega_k)^{n+1}=\mathbbm{1}_{D_k}(\tilde{\theta}^{n+1}+(n+1)\tilde{\theta}^n\wedge \omega_k)$$
This implies that for any $\epsilon>0$: 
\begin{eqnarray*}
(n+1)\int_{\mathbb{P}^1}\left(\int_{X_{\tau}\cap D_k}\tilde{\theta}^n\right)\omega_k\geq \int_{\{|\tau|\leq \epsilon\}}MA_{\tilde{\theta}}(u_k)-\int_{\{|\tau|\leq \epsilon\}}\mathbbm{1}_{D_k}\tilde{\theta}^{n+1}=\\=\int_{\{|\tau|\leq \epsilon\}}MA_{\tilde{\theta}}(u_k)-o(\epsilon).
\end{eqnarray*}
Combined with Proposition \ref{propliminf} we get the result.
\end{proof}

\begin{corollary} \label{corkey}
There exists a sequence $\tau_k$ such that $e^{k/2}\tau_k\to 0$ and $$\int_{X_{\tau_k}\cap D_k}\tilde{\theta}^n\to \textrm{vol}(\alpha).$$
\end{corollary}

\begin{proof}
We note from above that $$\omega_k=(e^{-k/2})_*\omega.$$ Also, since $\chi''>0$ for $x<-1$, $\int_{|\tau|<\delta}\omega>0$ for all $\delta>0$. Also note that for each $\tau$: $$\int_{X_{\tau_k}\cap D_k}\tilde{\theta}^n\leq \int_{X_{\tau}}MA_{\tilde{\theta}_{|X_{\tau}}}({u_k}_{|X_{\tau}})= \textrm{vol}(\alpha).$$That we can find a sequence $\tau_k$ with the desired properties now clearly follows from Proposition \ref{proplimit}.
\end{proof}

\section{Step 4: Convergence} \label{Sec:conv}

Let $\tilde{u}$ be as in Section \ref{Sec:form}, let $\eta$ be a K\"ahler form on $\mathcal{X}$ and pick an $\epsilon$. By Demailly's regularization theorem there exists an $\tilde{\theta}+\epsilon\eta$-psh function $\psi$ with analytic singularities such that $\psi\geq \tilde{u}$. We can assume that $E_T\subseteq E_{nK}(\beta)$. 

Let $\pi: \mathcal{Z}\to \mathcal{X}$ be a smooth modification with center contained in $E_{nK}(\beta)$ so that $\pi\circ \psi$ has divisorial singularities. It means that $$dd^c(\pi \circ \psi)=\sum_i a_i[E_i]+\sigma$$ where $\sigma$ is a smooth form. Since $$dd^c(\pi \circ \psi)\geq -\pi^*(\tilde{\theta}+\epsilon\eta)$$ we get that $\pi^*(\tilde{\theta}+\epsilon\eta)+\sigma$ is semipositive. For small $\delta>0$ the class $\pi^*\gamma-\delta\sum_i\{E_i\}$ is K\"ahler, so by possibly increasing $\epsilon$ by an arbitrarily small amount while adding arbitrarily small multiples of $\ln|s_{E_i}|^2_{h_i}$ to $\psi$ one can make sure that $$\tilde{\omega}:=\pi^*(\tilde{\theta}+\epsilon\eta)+\sigma$$ is K\"ahler.

Let us now consider two families of envelopes: $$v_k:=\sup\{v\leq g_k, \phi\preceq \psi: \phi\in PSH(\mathcal{X},\tilde{\theta}+\epsilon\eta) \},$$ and $$\tilde{v}_k:=\sup\{\tilde{v}\leq g_k\circ \pi-\psi\circ \pi: \tilde{\phi}\in PSH(\mathcal{Z},\tilde{\omega}) \}.$$ The two families are related in a simple way, namely:

\begin{proposition} \label{blabla}
$$v_k\circ \pi=\tilde{v}_k+\psi\circ \pi.$$
\end{proposition}

\begin{proof}
This is clear since a function $v$ is a candidate for the supremum $v_k$ iff $\tilde{v}:=v \circ \pi-\psi\circ\pi$ is a candidate for the supremum $\tilde{v}_k$.
\end{proof}

\begin{proposition}
${v_k}_{|X'}+k$ increases almost everywhere to a function $v_{X'}\in PSH(X',(\tilde{\theta}+\epsilon\eta)_{|X'})$ and similarly ${v_k}_{|\mathcal{E}}+k$ increases to a function $v_{\mathcal{E}}\in PSH(\mathcal{E},(\tilde{\theta}+\epsilon\eta)_{|\mathcal{E}})$. 
\end{proposition}

\begin{proof}
This follows since $g_k+k=\chi(\ln|\tau|^2+k)-\ln(1+|\tau|^2)$ is increasing and hence $v_k+k$ is also increasing, while $g_k+k=-1$ on $X'$ and $\mathcal{E}$, meaning that ${v_k}_{|X'}+k$ and ${v_k}_{|\mathcal{E}}+k$ have bounded limits.
\end{proof}

\begin{proposition} \label{propimp}
$$\int_{X'}MA_{(\tilde{\theta}+\epsilon\eta)_{|X'}}(v_{X'})\leq \textrm{vol}_{\mathcal{X}|X'}(\beta+\epsilon[\eta]),$$ and $$\int_{\mathcal{E}}MA_{(\tilde{\theta}+\epsilon\eta)_{|\mathcal{E}}}(v_{\mathcal{E}})\leq \textrm{vol}_{\mathcal{X}|\mathcal{E}}(\beta+\epsilon[\eta]).$$
\end{proposition}

\begin{proof}
Note that by Proposition \ref{blabla} $|v_k-\psi|\leq C_k$ for some constants $C_k$, and since $v_k+k$ was increasing, we see that $$v_{X'}\geq \psi_{|X'}-C_0$$ and similarly $$v_{\mathcal{E}}\geq \psi_{|\mathcal{E}}-C_0.$$  By the monotonicity of the Monge-Amp\`ere we thus get that $$\int_{X'}MA_{(\tilde{\theta}+\epsilon\eta)_{|X'}}(\psi_{|X'})\leq\int_{X'}MA_{(\tilde{\theta}+\epsilon\eta)_{|X'}}(v_{X'})$$ and $$\int_{\mathcal{E}}MA_{(\tilde{\theta}+\epsilon\eta)_{|\mathcal{E}}}(\psi_{|\mathcal{E}})\leq\int_{\mathcal{E}}MA_{(\tilde{\theta}+\epsilon\eta)_{|\mathcal{E}}}(v_{\mathcal{E}}).$$  

Since ${v_k}_{|X'}+k$ increases almost everywhere to $v_{X'}$ it follows that $MA_{(\tilde{\theta}+\epsilon\eta)_{|X'}}(v_k)$ converges weakly to $MA_{(\tilde{\theta}+\epsilon\eta)_{|X'}}(v_{X'})$ on the open sets $X'\cap\{\psi>-C\}$. It follows that $$\liminf_{k\to\infty}\int_{X'\cap \{\psi>-C\}}MA_{(\tilde{\theta}+\epsilon\eta)_{|X'}}({v_k}_{|X'})\geq \int_{X'\cap \{\psi>-C\}}MA_{(\tilde{\theta}+\epsilon\eta)_{|X'}}(v_{X'}).$$ Letting $C\to\infty$ we get that $$\liminf_{k\to\infty}\int_{X'}MA_{(\tilde{\theta}+\epsilon\eta)_{|X'}}({v_k}_{|X'})\geq \int_{X'}MA_{(\tilde{\theta}+\epsilon\eta)_{|X'}}(v_{X'}).$$ But since  $|v_k-\psi|\leq C_k$ we have that for all $k$ $$\int_{X'}MA_{(\tilde{\theta}+\epsilon\eta)_{|X'}}({v_k}_{|X'})=\int_{X'}MA_{(\tilde{\theta}+\epsilon\eta)_{|X'}}(\psi_{|X'}),$$ which shows that $$\int_{X'}MA_{(\tilde{\theta}+\epsilon\eta)_{|X'}}(v_{X'})=\int_{X'}MA_{(\tilde{\theta}+\epsilon\eta)_{|X'}}(\psi_{|X'}).$$ Since $$\int_{X'}MA_{(\tilde{\theta}+\epsilon\eta)_{|X'}}(\psi_{|X'})\leq \textrm{vol}_{\mathcal{X}|X'}(\beta+\epsilon[\eta])$$ we get that $$\int_{X'}MA_{(\tilde{\theta}+\epsilon\eta)_{|X'}}(v_{X'})\leq \textrm{vol}_{\mathcal{X}|X'}(\beta+\epsilon[\eta]),$$ and the same argument works for $\mathcal{E}$.
\end{proof}

The goal is now to prove the following.

\begin{proposition} \label{Prop:conv}
Given a sequence $\tau_k$ such that $e^{k/2}\tau_k\to 0$, then for any subsequence $\tau_{k_l}$ such that $e^{k_l/2}\tau_{k_l}$ is decreasing and $$\sum_{l=1}^{\infty}le^{k_l/2}\tau_{k_l}<\infty,$$ we have that the measures $MA_{(\tilde{\theta}+\epsilon\eta)_{|X_{\tau_{k_l}}}}({v_{k_l}}_{|X_{\tau_{k_l}}})$ converge locally to $MA_{(\tilde{\theta}+\epsilon\eta)_{|X'}}(v_{X'})$ and $MA_{(\tilde{\theta}+\epsilon\eta)|_{|\mathcal{E}}}(v_{|\mathcal{E}})$ away from $X'\cap \mathcal{E}$ and the singularities of $\psi$.
\end{proposition}

To prove this proposition we will use Berman's regularity theorem (Theorem \ref{ThmBerman}) to derive a Lipschitz bound for $v_k$. 

Let $x\in X\setminus Y$ and pick local holomorphic coordinates $z_i$ centered at $x$ defined in some neighbourhood $U$ of $x$. We then have that $(z_1,...,z_n,\tau)$ are coordinates for $U\times \mathbb{C}\subseteq \mathcal{X}$, let $U_k:=U\times \{|\tau|\leq e^{-k/2}\}.$ Let us also assume that $\psi$ is bounded on $U_0$.

In what follows $C$ will stand for a constant that can depend on the particular data, i.e. $\theta,\eta,\psi, U$ and so on, but importantly it will not depend on $k$.

\begin{lemma} \label{keyreglemma}
On $U_k$ we have that $$|v_k+k|\leq C$$ and $$|\Delta_{\tilde{\omega}}v_k|\leq Ce^k.$$
\end{lemma}

\begin{proof}
We have that $\psi-\max \psi-1\leq g$ and hence $\psi-\max\psi-1\leq v_0$ and since $\psi$ was assumed to be bounded on $U_0$ we get that $-C\leq v_0$ on $U_0$. Furthermore, since $v_k+k$ is increasing in $k$ we get that $-C\leq v_k+k$ on $U_k\subseteq U_0$ for all $k$. We also note that $v_k+k\leq g_k+k\leq 0$ on $U_k$, so we get that $$|v_k+k|\leq C$$ on $U_k$.

From Theorem \ref{ThmBerman} we get that $$\Delta_{\tilde{\omega}}\tilde{v}_k\leq (C(n+1)+\sup_{\mathcal{Z}}(\Delta_{\tilde{\omega}}g_k\circ \pi))e^{B(g_k\circ \pi-\psi\circ \pi-\inf_{\mathcal{Z}}(g_k\circ \pi-\psi\circ \pi))}.$$ Since $g_k(\tau)=\chi(\ln|e^{k/2}\tau|^2)-k-\ln(1+|\tau|^2)$ it follows that $$\sup_{\mathcal{Z}}(\Delta_{\tilde{\omega}}g_k\circ \pi)\leq Ce^k.$$ Since $v_k\circ \pi=\tilde{v}_k+\psi\circ \pi$ and $\pi$ is biholomorphic on $\pi^{-1}(U_0)$ ($\psi$ being smooth there) we get that on $U_0$: $$\Delta_{\tilde{\omega}} v_k\leq Ce^ke^{B(g_k-\psi-\inf(g_k-\psi))}.$$ Using that $\psi$ is bounded on $U_0$, $g_k-\inf g_k=g_k+k+1$ which is less than one on $U_k$ it follows that on $U_k$: $$\Delta_{\tilde{\omega}} v_k\leq Ce^k.$$ That $|\Delta_{\tilde{\omega}} v_k|\leq Ce^k$ is then immediate since $v_k$ is $\tilde{\theta}+\epsilon\eta$-psh.
\end{proof}

\begin{corollary} \label{corestimate}
On $U\times \{|\tau|\leq \frac{e^{-k/2}}{2}\}$ we have that
$$|v_k(z,\tau)-v_k(z,0)|\leq Ce^{k/2}|\tau|.$$ 
\end{corollary}

\begin{proof}
For a given $z\in U$ let $f_k(w):=v_k(z,e^{-k/2}w)+k$. Note that on $U_0$, $\tilde{\omega}$ is equivalent to the standard metric on $U\times \mathbb{D}$. It thus follows from Lemma \ref{keyreglemma} that $$|f_k|\leq C$$ and $$|\Delta f_k|\leq e^{-k}C|\Delta_{\tilde{\omega}}v_k|\leq C.$$ It then follows from Riesz representation theorem that for $|w|\leq 1/2$: $$|f_k(w)-f_k(0)|\leq C|w|.$$
\end{proof}

Let now $\tilde{z}_1,z_2,...,z_n,\tau$ be local coordinates centered at a point $z\in \mathcal{E}$, and assume that $\psi(z)>-\infty$. Let $V\subseteq \mathcal{E}$ be a neighbourhood of $z$ on which we assume that $\psi$ is bounded. The same arguments as for $X'$ shows that on $V\times \{|\tau|\leq \frac{e^{-k/2}}{2}\}$ (in local coordinates) we have that
$$|v_k(z,\tau)-v_k(z,0)|\leq Ce^{k/2}|\tau|.$$ 

We can now prove Proposition \ref{Prop:conv}.

\begin{proof}[Proof of Proposition \ref{Prop:conv}]
Using Corollary \ref{corestimate} (and the corresponding estimate near $\mathcal{E}$) and the fact that $v_k(z,0)+k$ is increasing in $k$ we get that if $k>m, |\tau_k|\leq e^{-k/2}/2$ and  $|\tau_m|\leq e^{-m/2}/2$ then 
\begin{eqnarray} \label{ineq1}
v_k(z,\tau_k)\geq v_k(z,0)-Ce^{k/2}|\tau_k|\geq v_m(z,0)-Ce^{k/2}|\tau_k|\geq \\ \geq v_m(z,\tau_m)-Ce^{k/2}|\tau_k|-Ce^{m/2}|\tau_m|.
\end{eqnarray}

This shows that $$v_{k_{l}}(z,\tau_{k_l})+\sum_{m=1}^l me^{k_m/2}|\tau_{k_m}|$$ is increasing for $l>2C$, and then it increases almost everywhere to $$v_{X'}(z)+\sum_{l=1}^{\infty}le^{k_l/2}|\tau_{k_l}|.$$ By the continuity of the Monge-Amp\`ere under increasing sequences the Proposition follows.
\end{proof}

Since $u_k\leq g_k,$ $u_k\in PSH(\tilde{\theta}+\epsilon\eta)$ and $u_k\preceq \psi$ it follows that $u_k\leq v_k$, and hence 
$D_k\subseteq D'_k$. Using Theorem \ref{DNT} we thus get the estimate 
\begin{eqnarray} \label{MAest}
MA_{(\tilde{\theta}+\epsilon\eta)_{|X_{\tau}}}({v_k}_{|X_{\tau}})\geq \mathbbm{1}_{X_{\tau}\cap D'_k}(\tilde{\theta}+\epsilon\eta)^n\geq \mathbbm{1}_{X_{\tau}\cap D_k}\tilde{\theta}^n.
\end{eqnarray}

We are now ready to prove the hard inequality of Theorem \ref{keyeq}.

\begin{proof}[Proof of Theorem \ref{keyeq}, hard inequality]

We know from Corollary \ref{corkey} that there is a sequence $\tau_k$ such that $e^{k/2}\tau_k\to 0$ and $$\int_{X_{\tau_k}\cap D_k}\tilde{\theta}^n\to \textrm{vol}(\alpha).$$ Let us pick a subsequence $\tau_{k_l}$ such that $e^{k_l/2}\tau_{k_l}$ is decreasing and $$\sum_{l=1}^{\infty}le^{k_l/2}\tau_{k_l}<\infty.$$ Of course 
\begin{equation} \label{cm}
\int_{X_{\tau_{k_l}}\cap D_{k_l}}\tilde{\theta}^n\to \textrm{vol}(\alpha).
\end{equation}

By Proposition \ref{Prop:conv} the measures $MA_{(\tilde{\theta}+\epsilon\eta)_{|X_{\tau_{k_l}}}}({v_{k_l}}_{|X_{\tau_{k_l}}})$ converge locally to $MA_{(\tilde{\theta}+\epsilon\eta)_{|X'}}(v_{X'})$ and $MA_{(\tilde{\theta}+\epsilon\eta)|_{|\mathcal{E}}}(v_{|\mathcal{E}})$ away from $X'\cap \mathcal{E}$ and the singularities of $\psi$. But for any $\delta>0$ we can find an open neighbourhood $U$ of the set of singularities $E\cup (X_0\cap E_{nK}(\beta))$ such that for all $|\tau|>0$: 
\begin{equation} \label{eqeq}
\int_{X_{\tau}\cap U}|\tilde{\theta}^n|<\delta.
\end{equation}
Hence the local convergence away from $U$ together with the estimates (\ref{MAest}), (\ref{cm}) and (\ref{eqeq}) implies that 
$$\int_{X'}MA_{(\tilde{\theta}+\epsilon\eta)_{|X'}}(v_{X'})+\int_{\mathcal{E}}MA_{(\tilde{\theta}+\epsilon\eta)|_{|\mathcal{E}}}(v_{|\mathcal{E}})\geq \textrm{vol}(\alpha)-\delta,$$ and since $\delta>0$ was arbitrary 
$$\int_{X'}MA_{(\tilde{\theta}+\epsilon\eta)_{|X'}}(v_{X'})+\int_{\mathcal{E}}MA_{(\tilde{\theta}+\epsilon\eta)|_{|\mathcal{E}}}(v_{|\mathcal{E}})\geq \textrm{vol}(\alpha).$$  By Proposition \ref{propimp} this implies that 
\begin{equation} \label{lasteq}
\textrm{vol}_{\mathcal{X}|X'}(\beta+\epsilon[\eta])+\textrm{vol}_{\mathcal{X}|\mathcal{E}}(\beta+\epsilon[\eta])\geq \textrm{vol}(\alpha).
\end{equation}

Since $K(\beta)$ intersects $X'$ and $\mathcal{E}$ the restricted volumes along $X'$ and $\mathcal{E}$ are continuous at $\beta$. Hence we can let $\epsilon\to 0$ in (\ref{lasteq}) to get the desired inequality $$\textrm{vol}_{\mathcal{X}|X'}(\beta)+\textrm{vol}_{\mathcal{X}|\mathcal{E}}(\beta)\geq \textrm{vol}(\alpha).$$

\end{proof}

\section{K\"ahler deformations to normal bundles}

The goal of this Section is to prove Theorem A and B. The outline of the argument was given in the Introduction.

\begin{proof}[Proof of Theorem \ref{ThmA}]
Let $\alpha:=[\omega]$. Recall that $\mu':X'\to X$ was the blow-up of $X$ along $Y$, with exceptional divisor $E$. Let $\alpha':=\mu'^*\alpha$. Let $\delta(\alpha,Y)$ denote the pseudoeffective threshhold of $Y$, i.e. the supremum of $t$ such that $\alpha'-t\{E\}$ is pseudoeffective. By continuity the volume of $\alpha'-t\{E\}$ tends to zero as $t\to \delta(\alpha,Y)$. Thus, given $\epsilon>0$ we can pick a $c$ such that $$0<\textrm{vol}(\alpha'-c\{E\})<\epsilon\textrm{vol}(\alpha).$$ Also pick some $b>c$ and let $$\beta=(\pi_X\circ\mu)^*\alpha+b\{X_0\}-c\{\mathcal{E}\}.$$

That $K(\beta)$ intersects each fiber $X_{\tau}$, $\tau\neq 0$, as well as $X'$ and $\mathcal{E}$ follows exactly as Proposition \ref{Prop:kl}, and one also notes that $E_{nK}(\beta)\subseteq X'$. That $$\textrm{vol}_{\mathcal{X}|X_{\tau}}=\textrm{vol}(\alpha)=\int_X\omega^n$$ is also easily seen, e.g. by considering the current $$T:=(\pi_X\circ\mu)^*\omega+(b-c)(\pi_{\mathbb{P}^1}\circ \mu)^*\omega_{FS}+(c-\delta)[X'].$$ Since $\beta_{|X'}=\alpha'-c\{E\}$ we have that $$\textrm{vol}_{\mathcal{X}|X'}(\beta)\leq \textrm{vol}(\alpha'-c\{E\})\leq \epsilon\textrm{vol}(\alpha).$$ Since $\{X_{\tau}\}=\{X'\}+\{\mathcal{E}\}$ Corollary \ref{Maincor} gives us that $$\textrm{vol}_{\mathcal{X}|\mathcal{E}}(\beta)=\textrm{vol}(\alpha)-\textrm{vol}(\alpha-c\{Y\})\geq (1-\epsilon)\textrm{vol}(\alpha).$$ 

By the definition of restricted volume this means that one can find a K\"ahler current $T\in \beta$ with analytic singularities such that $$\int_{N_{Y|X}\setminus E_T}T^n>(1-2\epsilon)\textrm{vol}(\alpha).$$ Without loss of generality we can assume that $E_T=E_{nK}(\beta)\subseteq X'$. We can now take a smooth modification $\pi':\mathcal{X}'\to \mathcal{X}$ with center contained in $X'$ such that $\pi'^*T$ has divisorial singularities, i.e. $\pi'^*T=\Omega'+\sum a_i[E_i]$. Thus $\Omega'$ is a semipositive form on $\mathcal{X}'$ such that $[\Omega'_{|X_1}]=\alpha$ and $\int_{N_{Y|X}}\Omega'>(1-2\epsilon)\textrm{vol}(\alpha)$. Letting $\sigma$ be a K\"ahler form on $\mathcal{X}'$ such that $[\sigma_{|X_1}]=\alpha$ and taking $(1-\delta)\Omega'+\delta \sigma$ shows that without loss of generality we can assume $\Omega'$ to be K\"ahler, and by taking an average with respect to the circle action we can make $\Omega'$ $S^1$-invariant. 

Let us now write $\Omega'=dd^cf+\pi_X^*\omega$ on $X\times \{|\tau|>1/2\}$ where $f$ is smooth and $S^1$-invariant. Let $\max_{reg}$ be a regularized max-function. If the constant $C$ is chosen appropriately we get that $\max_{reg}(-f+\ln(1+|\tau|^2)+C,0)$ is equal to $0$ for $|\tau|\leq 1$. It follows that $$\Omega:=\Omega'+dd^c\max_{reg}(-f+\ln(1+|\tau|^2)+C,0)$$ is a well-defined $S^1$-invariant K\"ahler form on $\mathcal{X}$, with $\int_{N_{Y|X}}\Omega>(1-2\epsilon)\textrm{vol}(\alpha)$ and $\Omega_{|X_{\tau}}=\omega$ for $|\tau|$ large. By scaling of $\tau$ we can then make sure that $\Omega_{|X_1}=\omega$, and we are done. 
\end{proof}

\begin{proof}[Proof of Theorem B]
We now let $c=\delta(\alpha,Y)$, $b>c$ and let $$\beta=(\pi_X\circ\mu)^*\alpha+b\{X_0\}-c\{\mathcal{E}\}.$$ We still have that $E_{nK}(\beta)\subseteq X'$.

Pick a closed positive current $\Omega'\in \beta$ with minimal singularities. As $K(\beta)$ intersects $\mathcal{E}$ the restricted volume along $\mathcal{E}$ is continuous at $\beta$, so it is the limit of restricted volumes of $$(\pi_X\circ\mu)^*\alpha+b\{X_0\}-c\{\mathcal{E}\}$$ as $c\to \delta(\alpha,Y)$ from below. As we saw in the proof of Theorem A this limit is exactly $\int_X\omega^n$, thus $$\int_{N_{Y|X}}\Omega'^n=\int_X\omega^n.$$

To then find an $\Omega$ such that $$\int_{N_{X|Y}}\Omega^n=\int_X\omega^n$$ which in addition is $S^1$-invariant and with $\Omega_{|X_1}=\omega$ we argue exactly as in the proof of Theorem A.
\end{proof}

\section{A canonical weak K\"ahler deformation} \label{Sec:can}

We will now describe how to construct a canonical weak K\"ahler deformation $(\mathcal{X}_{\overline{\mathbb{D}}},\Omega_{can})$ of $(X,\omega)$ to the normal bundle of $Y$.

Let $$u:=\sup\{v\leq 0: v\in PSH(X\times \overline{\mathbb{D}},\pi_X^*\omega), \nu_{Y\times \{0\}}(v)\geq \delta(\alpha,Y)\}.$$ Then $u$ is $S^1$-invariant, $\pi^*\omega$-psh, and has Lelong number along $Y\times \{0\}$ at least $\delta(\alpha,Y)$. Since $\delta(\alpha,Y)\ln|\tau|^2$ is a candidate for the supremum, and $u\leq 0$, we get that $u_{X_1}=0$.  We now let $$\Omega_{can}:=dd^c(u\circ \mu)+(\pi_X\circ \mu)^*\omega-\delta(\alpha,Y)[\mathcal{E}].$$ This is an $S^1$-invariant closed positive current, and since $u_{|X_1}=0$ we have that ${\Omega_{can}}_{|X_1}=\omega$. Thus $(\mathcal{X}_{\overline{\mathbb{D}}},\Omega_{can})$ is a weak K\"ahler deformation of $(X,\omega)$ to the normal bundle of $Y$.

Note that $$\delta(\alpha,Y)\ln|\tau|^2\leq u\leq 0$$ implies that $\max(u,(\delta(\alpha,Y)+1)\ln|\tau|^2+1)$ is equal to $u$ for $\ln|\tau|^2\leq -1$ and $(\delta(\alpha,Y)+1)\ln|\tau|^2+1$ when $(\delta(\alpha,Y)+1)\ln|\tau|^2+1\geq 0$. It follows that $\max(u,(\delta(\alpha,Y)+1)\ln|\tau|^2+1)$ extends as a $\pi^*\omega$-psh function on $X\times \mathbb{C}$, which is equal to $u$ for $|\tau|$ small.

Let as in Section \ref{Sec:form} $s_{\mathcal{E}}$ be a defining section for $\mathcal{E}$, $h$ a smooth metric, and $\eta$ the form so that $$dd^c\ln|s_{\mathcal{E}}|^2_h=[\mathcal{E}]-\eta.$$ Let $$\tilde{\theta}:=(\pi_X\circ \mu)^*\omega+(\delta(\alpha,Y)+1)(\pi_{\mathbb{P}^1}\circ \mu)^*\omega_{FS}-\delta(\alpha,Y)\eta.$$

It follows that $$\tilde{u}:=\max(u,(\delta(\alpha,Y)+1)\ln|\tau|^2+1)-(\delta(\alpha,Y)+1)\ln(1+|\tau|^2)-\delta(\alpha,Y)\ln|s_{\mathcal{E}}|^2_h$$ is $\tilde{\theta}$-psh and that $dd^c\tilde{u}+\tilde{\theta}=\Omega_{can}$ for $\ln|\tau|^2\leq -1$. 

We now claim that $\tilde{u}$ has minimal singularities. To see this, let $\tilde{v}$ be an $\tilde{\theta}$-psh function. Then $$v:=\tilde{v}+(\delta(\alpha,Y)+1)\ln(1+|\tau|^2)+\delta(\alpha,Y)\ln|s_{\mathcal{E}}|^2_h$$ is $(\pi_X\circ \mu)^*\omega$-psh and has Lelong number at least $\delta(\alpha,Y)$ along $\mathcal{E}$, so it descends to a $\pi_X^*\omega$-psh function on $X\times \mathbb{C}$ with Lelong number at least $\delta(\alpha,Y)$ along $Y\times \{0\}$. Since it is also automatically bounded from above by some constant, we get that $v\leq u+C$ for $|\tau|\leq 1$. For $|\tau|\geq 1$: $\tilde{u}$ is bounded, hence $\tilde{v}\leq \tilde{u}+C'$, showing that $\tilde{u}$ has minimal singularities.

It follows that $$\int_{N_{X|Y}}\Omega_{can}^n=\textrm{vol}_{\mathcal{X}|\mathcal{E}}(\beta),$$ where $$\beta:=(\pi_X\circ \mu)^*\alpha+(\delta(\alpha,Y)+1)\{X_0\}-\delta(\alpha,Y)\{\mathcal{E}\}.$$ But we saw that $$\textrm{vol}_{\mathcal{X}|\mathcal{E}}(\beta)=\int_X\omega^n$$ already in the proof of Theorem \ref{ThmB}.

\section{Proof of Theorem A' and B'}

Let $\gamma$ be a K\"ahler class on $\mathcal{Z}$ such that $\gamma_{X_1}=\alpha,$ and let $\eta$ be a K\"ahler form in $\gamma.$ Let $\beta:=\gamma+b\{X_{\infty}\}-c\{Z_0\}$ where $b>c>0$. As in the proof of Theorem A we see that $E_{nK}(\beta)\subseteq D$. We claim that for $c$ large enough $\nu_{Z_i}(\beta)>0$ for all $i>0$. It would then follow from Corollary A that $\textrm{vol}_{\mathcal{Z}|Z_0}(\beta)=\int_X\omega^n,$ and then the proofs of Theorem A' and B' follows exactly those of Theorem A and B'.

\begin{proof}[Proof of claim]
Let $\delta>0$. We can then find a closed positive current $T\in \beta+\delta\gamma$ with analytic singularities such that $\nu_{Z_0}(T)=0$  and $a_i:=\nu_{Z_i}(T)\leq \nu_{Z_i}(\beta)$ and  for all $i>0$. Let $a_0:=c$. 

Fix an $i$. Let $Z_{i,j}:=Z_i\cap Z_j$ and let $A$ be the set indices $j$ different from $i$ such that $Z_{i,j}\neq \emptyset.$ Let also $B$ be the subset of $A$ such that $a_j\geq a_i$. Let $$T_i:=(T+c[Z_0]-a_i[\mathcal{Z}_0])_{|Z_i},$$ where $\mathcal{Z}_0$ denotes the zero-fiber of $\mathcal{Z}$. Note that $[T_i]=(1+\delta)\gamma_{|Z_i},$ thus $$\int_{Z_i}T_i\wedge \eta^{n-1}=(1+\delta)\int_{Z_i}\eta^n.$$ On the other hand $T_i=\sigma+\sum_{j\in A}(a_{j}-a_i)[Z_{i,j}]$ where $\sigma$ is semipositive. Hence we get the inequality 
\begin{eqnarray} \label{inin}
(1+\delta)\int_{Z_i} \eta^{n} \geq \sum_{j\in A}(a_j-a_i)\int_{Z_{i,j}}\eta^{n-1}\geq   \nonumber \\ \geq \sum_{j\in B}(a_j-a_i)\int_{Z_{i,j}}\eta^{n-1}-a_i\sum_{j\in A\setminus B}\int_{Z_{i,j}}\eta^{n-1}.
\end{eqnarray}
If we pick one $l\in B$ and rearrange we get $$a_l\leq \frac{1}{\int_{Z_{i,l}}\eta^{n-1}}\left((1+\delta)\int_{Z_i}\eta^n+a_i\sum_{j\in (A\setminus B)\cap\{l\}}\int_{Z_{i,j}}\eta^{n-1}\right).$$ This shows that we can find constants $C_1,C_2$ independent of $\delta$ and $i$ such that for all $j\in A$ $$a_j\leq C_1a_i+C_2.$$ Since $\mathcal{Z}_0$ is connected this implies that if $a_0=c$ is large enough each $a_i>0$, and since this estimate was uniform in $\delta,$ each $\nu_{Z_i}(\beta)>0$ for $i>0$.
\end{proof}

\bigskip
\noindent David Witt Nystr\"om \newline
Department of Mathematical Sciences \newline
Chalmers University of Technology and the University of Gothenburg \newline
SE-412 96 Gothenburg, Sweden \newline
wittnyst@chalmers.se

\end{document}